# STABLE LIMITS OF MARTINGALE TRANSFORMS WITH APPLICATION TO THE ESTIMATION OF GARCH PARAMETERS

By Thomas Mikosch[1] and Daniel Straumann

*University of Copenhagen and ETH Zürich*

In this paper we study the asymptotic behavior of the Gaussian quasi maximum likelihood estimator of a stationary GARCH process with heavy-tailed innovations. This means that the innovations are regularly varying with index $\alpha \in (2, 4)$. Then, in particular, the marginal distribution of the GARCH process has infinite fourth moment and standard asymptotic theory with normal limits and $\sqrt{n}$-rates breaks down. This was recently observed by Hall and Yao [*Econometrica* **71** (2003) 285–317]. It is the aim of this paper to indicate that the limit theory for the parameter estimators in the heavy-tailed case nevertheless very much parallels the normal asymptotic theory. In the light-tailed case, the limit theory is based on the CLT for stationary ergodic finite variance martingale difference sequences. In the heavy-tailed case such a general result does not exist, but an analogous result with infinite variance stable limits can be shown to hold under certain mixing conditions which are satisfied for GARCH processes. It is the aim of the paper to give a general structural result for infinite variance limits which can also be applied in situations more general than GARCH.

**1. Introduction.** The motivation for writing this paper comes from Gaussian quasi maximum likelihood estimation (QMLE) for GARCH (*generalized autoregressive conditionally heteroscedastic*) processes with regularly varying noise; we refer to Section 4 for a detailed description of the problem.

Received November 2003; revised March 2005.

[1]Supported in part by MaPhySto, the Danish Research Network for Mathematical Physics and Stochastics, DYNSTOCH, a research training network under the Improving Human Potential Programme financed by the Fifth Framework Programme of the European Commission, and by Danish Natural Science Research Council (SNF) Grant 21-01-0546.

Recall that the process

$$(1.1) \qquad X_t = \sigma_t Z_t \qquad \text{with } \sigma_t^2 = \alpha_0 + \sum_{i=1}^{p} \alpha_i X_{t-i}^2 + \sum_{j=1}^{q} \beta_j \sigma_{t-j}^2, t \in \mathbb{Z},$$

is said to be a GARCH$(p, q)$ process [GARCH process of *order* $(p, q)$]. Here $(Z_t)$ is an i.i.d. sequence with $EZ_1^2 = 1$ and $EZ_1 = 0$, and $\alpha_i, \beta_j$ are non-negative constants. GARCH processes and their parameter estimation have been intensively investigated over the last few years; see [19] for a general overview and [28] and the references therein for parameter estimation in GARCH and related models. In the context of QMLE, the asymptotic behavior of the parameter estimator is essentially determined by the limiting behavior of the quantity [see (4.13)]

$$L_n'(\boldsymbol{\theta}_0) = \frac{1}{2} \sum_{t=1}^{n} \frac{h_t'(\boldsymbol{\theta}_0)}{\sigma_t^2} (Z_t^2 - 1),$$

where $L_n'$ is the derivative of the underlying log-likelihood, $h_t'$ is the derivative of $\sigma_t^2$ when considered as a function of the parameter $\boldsymbol{\theta}$, and $\boldsymbol{\theta}_0$ is the true parameter (consisting of the $\alpha_i$ and $\beta_j$ values) in a certain parameter space. In this context,

$$\mathbf{G}_t = \frac{h_t'(\boldsymbol{\theta}_0)}{\sigma_t^2}, \qquad t \in \mathbb{Z},$$

is a stationary ergodic sequence of vector-valued random variables which is adapted to the filtration $\mathcal{F}_t = \sigma(Y_{t-1}, Y_{t-2}, \dots)$, $t \in \mathbb{Z}$, where $Y_t = Z_t^2 - 1$ constitutes an i.i.d. sequence.

If $\mathbf{G}_t$ has a finite first moment, the sequence $(\mathbf{G}_t Y_t)$ is a transform of the martingale difference sequence $(Y_t)$, hence, a stationary ergodic martingale difference sequence with respect to $(\mathcal{F}_t)$. If $E|\mathbf{G}_1|^2 < \infty$ and $EY_1^2 < \infty$, an application of the central limit theorem (CLT) for finite variance stationary ergodic martingale differences (see [4], Theorem 23.1) yields

$$n^{-1/2} \sum_{t=1}^{n} \mathbf{G}_t Y_t \xrightarrow{d} N(\mathbf{0}, \boldsymbol{\Sigma}),$$

where $\boldsymbol{\Sigma}$ is the covariance matrix of $\mathbf{G}_1 Y_1$. This result does not require any additional information about the dependence structure of $(\mathbf{G}_t Y_t)$. It implies the asymptotic normality of the parameter estimator based on QMLE.

If $EY_1^2 = \infty$, a result as general as the CLT for stationary ergodic martingale differences is not known. However, some limit results for stationary sequences with marginal distribution in the domain of attraction of an infinite variance stable distribution exist. We recall two of them in Section 2. Our interest in infinite variance stable limit distributions for $\sum_{t=1}^{n} \mathbf{G}_t Y_t$ is again closely related to parameter estimation for GARCH processes. Recently, Hall and Yao [16] gave the asymptotic theory for QMLE in GARCH models when $EZ_1^4 = \infty$. To be more specific, they assume regular variation with index $\alpha \in (1, 2)$ for the distribution of $Z_1^2$. It is our aim to show that their results can be obtained by a general limit result for the martingale transforms $\sum_{t=1}^{n} \mathbf{G}_t Y_t$ when the i.i.d. noise $(Y_t)$ is regularly varying with index $\alpha \in (1, 2)$. The key notions in this context are *regular variation* of the finite-dimensional distributions of $(\mathbf{G}_t Y_t)$ and strong mixing of this sequence; see Section 2 for these notions.

Our objective is twofold. *First*, we want to show that the theories on parameter estimation for GARCH processes with heavy- or light-tailed innovations $(Z_t)$ parallel each other. We use the recent structural approach to GARCH estimation by Berkes et al. [3] in order to show that such a unified approach is possible. *Second*, our approach to the asymptotic theory for parameter estimators is not restricted to GARCH processes. In the light-tailed case, Straumann and Mikosch [28] extended the approach by Berkes et al. [3], including among others AGARCH and EGARCH processes. The main difficulty of our approach when infinite variance limits occur is the verification of certain mixing conditions. In contrast to the case of asymptotic normality, such conditions cannot be avoided. However, it is difficult to check for a given model that these conditions hold; see Section 4.4 in order to get a flavor of the task to be solved.

GARCH processes and their parameter estimation give the motivation for this paper. The corresponding limit theory for the QMLE with heavy-tailed innovations can be found in Section 4. Our main tool for achieving these limit results is based on asymptotic theory for martingale transforms with infinite variance stable limits. This theory is formulated and proved in Section 3. It is based on more general results for sums of stationary mixing vector sequences with regularly varying finite-dimensional distributions. This theory is outlined in Section 2.

**2. Preliminaries.** In this section we collect some basic tools and notions to be used throughout this paper. First we want to formulate a classical result on infinite variance stable limits for i.i.d. vector-valued summands due to Rvačeva [25]. Before we formulate this result, we recall the notions of *stable random vector* and *multivariate regular variation*. The class of stable random vectors coincides with the class of possible limit distributions for sums of i.i.d. random vectors, and multivariate regular variation is the domain of attraction condition for sums of i.i.d. random vectors. Then we continue with an analog of Rvačeva's result for stationary ergodic vector sequences. In this context, we also need to recall some *mixing conditions*.

*Stable random vectors.*  Recall that a vector $\mathbf{X}$ with values in $\mathbb{R}^d$ is said to be $\alpha$-stable for some $\alpha \in (0, 2)$ if its characteristic function is given by

$$
(2.1)
\begin{aligned}
&E e^{i(\mathbf{x}, \mathbf{X})} \\
&= \begin{cases}
\exp\Big\{ -\int_{\mathbb{S}^{d-1}} |(\mathbf{x}, \mathbf{y})|^\alpha (1 - i \operatorname{sign}((\mathbf{x}, \mathbf{y})) \tan(\pi\alpha/2)) \\
\qquad\qquad\qquad\qquad\qquad \times \Gamma(d\mathbf{y}) + i(\mathbf{x}, \boldsymbol{\mu}) \Big\}, \qquad \alpha \neq 1, \\
\exp\Big\{ -\int_{\mathbb{S}^{d-1}} |(\mathbf{x}, \mathbf{y})| \Big( 1 + i\frac{2}{\pi} \operatorname{sign}((\mathbf{x}, \mathbf{y})) \log |(\mathbf{x}, \mathbf{y})| \Big) \\
\qquad\qquad\qquad\qquad\qquad \times \Gamma(d\mathbf{y}) + i(\mathbf{x}, \boldsymbol{\mu}) \Big\}, \qquad \alpha = 1,
\end{cases}
\end{aligned}
$$

where $(\mathbf{x}, \mathbf{y})$ denotes the usual inner product in $\mathbb{R}^d$ and $|\cdot|$ the Euclidean norm; see [27], Theorem 2.3.1. The *index of stability* $\alpha \in (0, 2)$, the *spectral measure* $\Gamma$ on the unit sphere $\mathbb{S}^{d-1}$ and the location parameter $\boldsymbol{\mu}$ uniquely determine the distribution of an infinite variance $\alpha$-stable random vector $\mathbf{X}$.

*Multivariate regular variation.*  If $\mathbf{X}$ is $\alpha$-stable for some $\alpha \in (0, 2)$, it is regularly varying with index $\alpha$. This means the following. The random vector $\mathbf{X}$ with values in $\mathbb{R}^d$ is *regularly varying with index* $\alpha \geq 0$ if there exists a random vector $\boldsymbol{\Theta}$ with values in the unit sphere $\mathbb{S}^{d-1}$ of $\mathbb{R}^d$ such that for any $t > 0$, as $x \to \infty$,

$$
(2.2) \qquad \frac{P(|\mathbf{X}| > tx, \widetilde{\mathbf{X}} \in \cdot)}{P(|\mathbf{X}| > x)} \overset{v}{\to} t^{-\alpha} P(\boldsymbol{\Theta} \in \cdot),
$$

where for any vector $\mathbf{x} \neq \mathbf{0}$,

$$
\tilde{\mathbf{x}} = \mathbf{x}/|\mathbf{x}|,
$$

and $\overset{v}{\to}$ denotes *vague convergence* in the Borel $\sigma$-field of $\mathbb{S}^{d-1}$; see [22, 23] for its definition and details. The distribution of $\boldsymbol{\Theta}$ is called the *spectral measure* of $\mathbf{X}$. Alternatively, (2.2) is equivalent to

$$
(2.3) \qquad \frac{P(\mathbf{X} \in x\cdot)}{P(|\mathbf{X}| > x)} \overset{v}{\to} \mu,
$$

where $\overset{v}{\to}$ denotes vague convergence in the Borel $\sigma$-field of $\overline{\mathbb{R}}^d \setminus \{\mathbf{0}\}$ and $\mu$ is a measure on the same $\sigma$-field satisfying the homogeneity assumption $\mu(tA) = t^{-\alpha} \mu(A)$ for $t > 0$.

REMARK 2.1.  The property of regular variation of $\mathbf{X}$ with index $\alpha$ does not depend on the chosen norm. However, the spectral measure (the unit spheres $\mathbb{S}^{d-1}$ depend on the norm) and the limiting measure $\mu$ can be different for distinct norms. The asymptotic theory of this paper does not depend on the particular choice of the norm $|\cdot|$. Unless specified otherwise, we will, however, assume that $|\cdot|$ is the Euclidean norm.

To give some intuition on regular variation of a vector $\mathbf{X}$, we mention some immediate consequences of the definition. Regular variation of $\mathbf{X}$ implies that $|\mathbf{X}|$ is regularly varying: $P(\mathbf{X}| > x) = L(x)x^{-\alpha}$, where $L(x)$ is slowly varying in the sense that $L(cx)/L(x) \to 1$ as $x \to \infty$, for every $c > 0$. This property follows by plugging the set $\mathbb{S}^{d-1}$ into (2.2). Moreover, relation (2.3) implies that every linear combination $(\mathbf{a}, \mathbf{X})$, $\mathbf{a} \neq \mathbf{0}$, of the components of $\mathbf{X}$ is regularly varying with the same index $\alpha$. This follows by plugging the $d$-dimensional halfspace $\{\mathbf{x} \in \mathbb{R}^d : (\mathbf{a}, \mathbf{x}) > 1\}$ into (2.3).

Definition (2.2) has an equivalent sequential analog in the following sense. Choosing any sequence $a_n \to \infty$ such that

$$(2.4) \qquad nP(|\mathbf{X}| > a_n) \to 1,$$

(2.2) is equivalent to

$$(2.5) \qquad nP(|\mathbf{X}| > ta_n, \widetilde{\mathbf{X}} \in S) \to t^{-\alpha} P(\mathbf{\Theta} \in S), \qquad t > 0,$$

for all Borel sets $S \subset \mathbb{S}^{d-1}$ with $P(\mathbf{\Theta} \in \partial S) = 0$. By an application of Poisson's limit theorem, the latter relation implies for an i.i.d. sequence $(\mathbf{X}_i)$ with the same marginal distribution as $\mathbf{X}$ that the binomial random variable

$$
\begin{aligned}
(2.6) \qquad & N_n((t, \infty) \times S) \\
& = \sum_{i=1}^{n} I_{(t,\infty) \times S}((a_n^{-1}|\mathbf{X}_i|, \widetilde{\mathbf{X}}_i)) \xrightarrow{d} N((t, \infty) \times S),
\end{aligned}
$$

where the limiting variable is Poisson with parameter $t^{-\alpha} P(\mathbf{\Theta} \in S)$ and $I_A$ denotes the indicator function of $A$. This binomial variable counts those exceedances of the scaled lengths $a_n^{-1}|\mathbf{X}_1|, \ldots, a_n^{-1}|\mathbf{X}_n|$ of the vectors $\mathbf{X}_i$ above the threshold $t$ for which the angles of the $\mathbf{X}_i$'s fall into the set $S$. The distributional convergence (2.6) can be extended to the weak convergence of the underlying point processes $N_n$ toward a Poisson process $N$ on $\overline{\mathbb{R}}^d \setminus \{\mathbf{0}\}$, $\mu$ being its mean measure; we omit the details and refer again to the mentioned literature [22, 23]. However, the limit relation (2.6) already explains to some extent what the spectral measure describes (in an asymptotic sense): it gives the likelihood that the angles of the i.i.d. regularly varying vectors $\mathbf{X}_1, \ldots, \mathbf{X}_n$ "far away from the origin" fall into a specified set $S$.

The Poisson convergence result (2.6) also tells us what "far away from the origin" means: the scaling $a_n$ of the $\mathbf{X}_i$'s has to be chosen according to the condition (2.4). We see in the sequel that this condition will appear in various disguises. Finally, we mention that (2.3) can be written in equivalent sequential form with $(a_n)$ satisfying (2.4) as

$$nP(a_n^{-1}\mathbf{X} \in \cdot) \xrightarrow{v} \mu(\cdot).$$

*Stable limits for sums of i.i.d. random vectors.*   Now let $(\mathbf{Y}_t)$ be an i.i.d. sequence of random vectors with values in $\mathbb{R}^d$. According to Rvačeva [25], there exist sequences of constants $a_n > 0$ and $\mathbf{b}_n \in \mathbb{R}^d$ such that

$$a_n^{-1} \sum_{t=1}^{n} \mathbf{Y}_t - \mathbf{b}_n \xrightarrow{d} \mathbf{X}_\alpha$$

for some $\alpha$-stable random variable $\mathbf{X}_\alpha$ with $\alpha \in (0,2)$ if and only if $\mathbf{Y}_1$ is regularly varying with index $\alpha$, and the normalizing constants $a_n$ can be chosen as

$$(2.7) \qquad\qquad P(|\mathbf{Y}_1| > a_n) \sim n^{-1}.$$

Notice that (2.7) is directly comparable with condition (2.4), which appears in the sequential definition of regular variation.

For a stationary sequence $(\mathbf{Y}_t)$, a similar result can be found in [13] as a multivariate extension of one-dimensional results in [12]. For its formulation one needs regular variation of the summands and a particular mixing condition, called $\mathcal{A}(a_n)$, which was introduced in [12].

*Mixing conditions.*   We say that the *condition $\mathcal{A}(a_n)$* holds for the stationary sequence $(\mathbf{Y}_t)$ of random vectors with values in $\mathbb{R}^d$ if there exists a sequence of positive integers $r_n$ such that $r_n \to \infty$, $k_n = [n/r_n] \to \infty$ as $n \to \infty$ and

$$(2.8) \quad E \exp\left\{ -\sum_{t=1}^{n} f(\mathbf{Y}_t/a_n) \right\} - \left( E \exp\left\{ -\sum_{t=1}^{r_n} f(\mathbf{Y}_t/a_n) \right\} \right)^{k_n} \to 0,$$

$$n \to \infty, \ \forall f \in \mathcal{G}_s,$$

where $\mathcal{G}_s$ is the collection of bounded nonnegative step functions on $\overline{\mathbb{R}}^d \setminus \{\mathbf{0}\}$. The convergence in (2.8) is not required to be uniform in $f$. This is indeed a very weak condition and is implied by many known mixing conditions, in particular, the strong mixing condition which is relevant in the context of GARCH processes; see Section 4. We refer to [13] for a comparison of $\mathcal{A}(a_n)$ with other mixing conditions.

For later use we also recall the definition of a *strongly mixing stationary sequence* $(\mathbf{Y}_t)$ of random vectors with *rate function* $(\phi_k)$ (see [24], cf. [14] or [17]):

$$\sup_{A \in \sigma(\mathbf{Y}_s, s \leq 0), B \in \sigma(\mathbf{Y}_s, s > k)} |P(A \cap B) - P(A)P(B)| =: \phi_k \to 0 \qquad \text{as } k \to \infty.$$

If $(\phi_k)$ decays to zero at an exponential rate, then $(\mathbf{Y}_t)$ is said to be *strongly mixing with geometric rate*. In Section 4.4 we use a more stringent notion of mixing, called $\beta$-mixing or absolute regularity. It implies strong mixing with the same rate function.

*Stable limits for sums of stationary random variables.* The following result is a combination of Theorem 2.8 and Proposition 3.3 in [13]. It gives conditions under which an $\alpha$-stable weak limit occurs for the sum process of a stationary sequence. In what follows we write

$$\mathbf{S}_0 = \mathbf{0} \quad \text{and} \quad \mathbf{S}_n = \mathbf{Y}_1 + \cdots + \mathbf{Y}_n, \qquad n \geq 1,$$

and for any Borel set $B \subset \mathbb{R}$,

$$\mathbf{S}_n B = (S_n^{(h)}(B))_{h=1,\ldots,d},$$

where

$$S_n^{(h)}(B) = \sum_{t=1}^n Y_t^{(h)} I_B(|Y_t^{(h)}|/a_n), \qquad n \geq 1.$$

THEOREM 2.2. *Let $(\mathbf{Y}_t)$ be a strictly stationary sequence of random vectors with values in $\mathbb{R}^d$ and the real sequence $(a_n)$ be defined by (2.7). Assume that the following conditions are satisfied:*

(a) *The finite-dimensional distributions of $(\mathbf{Y}_k)$ are regularly varying with index $\alpha > 0$. To be specific, let $\text{vec}(\boldsymbol{\theta}_{-k}^{(k)}, \ldots, \boldsymbol{\theta}_k^{(k)})$ be the $(2k+1)d$-dimensional random row vector with values in the unit sphere $\mathbb{S}^{(2k+1)d-1}$ that appears in the definition (2.2) of regular variation of $\text{vec}(\mathbf{Y}_{-k}, \ldots, \mathbf{Y}_k)$, $k \geq 0$, with respect to the max-norm $|\cdot|$ in $\mathbb{R}^{(2k+1)d}$.*

(b) *The mixing condition $\mathcal{A}(a_n)$ holds for $(\mathbf{Y}_t)$.*

(c)

$$\lim_{k \to \infty} \limsup_{n \to \infty} P\left( \bigvee_{k \leq |t| \leq r_n} |\mathbf{Y}_t| > a_n y \,\Big|\, |\mathbf{Y}_0| > a_n y \right) = 0, \qquad y > 0, \tag{2.9}$$

*where $(r_n)$ appears in the formulation of $\mathcal{A}(a_n)$.*

*Then the limit*

$$\gamma = \lim_{k \to \infty} E\left( |\boldsymbol{\theta}_0^{(k)}|^\alpha - \bigvee_{j=1}^k |\boldsymbol{\theta}_j^{(k)}|^\alpha \right)_+ \Big/ E|\boldsymbol{\theta}_0^{(k)}|^\alpha \tag{2.10}$$

*exists. If $\gamma > 0$, then the following results hold:*

(i) *If $\alpha \in (0, 1)$, then*

$$a_n^{-1} \mathbf{S}_n \overset{d}{\to} \mathbf{X}_\alpha,$$

*for some $\alpha$-stable random vector $\mathbf{X}_\alpha$.*

(ii) *If $\alpha \in [1, 2)$, and for all $\delta > 0$*

$$(2.11) \qquad \lim_{y \to 0} \limsup_{n \to \infty} P(|\mathbf{S}_n(0, y] - E\mathbf{S}_n(0, y]| > \delta a_n) = 0,$$

*then*

$$a_n^{-1}(\mathbf{S}_n - E\mathbf{S}_n(0, 1]) \xrightarrow{d} \mathbf{X}_\alpha$$

*for some $\alpha$-stable random vector $\mathbf{X}_\alpha$.*

REMARK 2.3.   The structure of the limiting vectors $\mathbf{X}_\alpha$ is given by some functional of the points of a limiting point process. The proof of this result makes heavy use of point process convergence results, which are appropriate tools in the context of regularly varying distributions when extremely large values may occur in the sequence $(\mathbf{Y}_t)$; see [13] for details. This leaves the parameters in the characteristic function (2.1) unspecified (with the exception of $\alpha$); a specification is not available so far and requires further investigation.

REMARK 2.4.   The quantity $\gamma$ in (2.4) can be identified as the *extremal index* of the sequence $(|\mathbf{Y}_t|)$; see [12] and Remark 2.3 in [13]. The extremal index $\gamma \in [0, 1]$ of a strictly stationary real-valued sequence is a number which characterizes the clustering behavior of the sequence above high thresholds. Roughly speaking, its existence ensures that the approximate relationship

$$P\left(\max_{i=1,\dots,n} |\mathbf{Y}_i| \le u_n\right) \approx P^{n\gamma}(|\mathbf{Y}_1| \le u_n)$$

holds for suitable sequences $u_n \to \infty$). For the definition and interpretation of the extremal index, we refer to [18] and [15], Section 8.1. The case $\gamma = 0$ corresponds to the case of sequences with unusually large cluster sizes above high thresholds. This case is often considered pathological; see [18] for some examples and the recent paper by Samorodnitsky [26]. For $\gamma = 0$ the limit theory developed in [12, 13] yields that the weak limit results in the above theorem hold with zero limit.

**3. Stable limits for martingale transform.**   In this section we want to derive infinite variance stable limits for sums of strictly stationary random vectors which have the particular form

$$\mathbf{Y}_t = \mathbf{G}_t Y_t,$$

where $(Y_t)$ is an i.i.d. sequence and $(\mathbf{G}_t)$ is a strictly stationary sequence of random vectors with values in $\mathbb{R}^d$ such that $(\mathbf{G}_t)$ is adapted to the filtration given by the $\sigma$-fields $\mathcal{F}_t = \sigma(Y_{t-1}, Y_{t-2}, \dots)$, $t \in \mathbb{Z}$. If $EY_1 = 0$ and $E|\mathbf{G}_1| <$ $\infty$, $E(\mathbf{G}_t Y_t | \mathcal{F}_t) = \mathbf{0}$ a.s., and, therefore, $(\mathbf{G}_t Y_t)$ is a martingale difference sequence and

$$\mathbf{S}_0 = \mathbf{0}, \qquad \mathbf{S}_n = \mathbf{Y}_1 + \cdots + \mathbf{Y}_n, \qquad n \geq 1,$$

is the *martingale transform* of the martingale $(\sum_{t=1}^n Y_t)_{n \geq 0}$ by the sequence $(\mathbf{G}_t)$. We keep this name even if $E|\mathbf{Y}_1| = \infty$.

3.1. *Basic assumptions.* We impose the following assumptions on the sequences $(Y_t)$ and $(\mathbf{G}_t)$:

A.1. $Y_1$ is regularly varying with index $\alpha \in (0, 2)$.
A.2. $E|\mathbf{G}_1|^{\alpha + \epsilon} < \infty$ for some $\epsilon > 0$.
A.3. $(\mathbf{G}_t Y_t)$ satisfies condition $\mathcal{A}(a_n)$ [see (2.8)], where $P(|Y_1| > a_n) \sim n^{-1}$ and $(r_n)$, defined in (2.8), is such that

$$(3.1) \qquad n r_n \left( \frac{a_{r_n}}{a_n} \right)^{\alpha + \epsilon} \to 0,$$

where $\epsilon$ is the same as in A.2.

REMARK 3.1. Regular variation of $Y_1$ with index $\alpha$ and the i.i.d. property of $(Y_t)$ imply that

$$P\left( a_n^{-1} \max_{1 \leq t \leq n} |Y_t| \leq x \right) \to \Phi_\alpha(x) = e^{-x^{-\alpha}}, \qquad x > 0,$$

for the Fréchet distribution $\Phi_\alpha$; see [15], Chapter 3.

In this setting, the heaviness of the tails of the distribution of $\mathbf{G}_1 Y_1$ is essentially determined by the distribution of $Y_1$; see Remark 3.4 below.

3.2. *Main result.* We are now ready to formulate our main result on the asymptotic behavior of the sum process $(\mathbf{S}_n)$.

THEOREM 3.2. *Consider the martingale transform*

$$\left( \sum_{t=1}^n \mathbf{Y}_t \right)_{n \geq 0} = \left( \sum_{t=1}^n \mathbf{G}_t Y_t \right)_{n \geq 0}$$

*defined above. Assume that the conditions* A.1–A.3 *are satisfied. Moreover, if $\alpha \in (1, 2)$, assume that $EY_1 = 0$ and, if $\alpha = 1$, that $Y_1$ is symmetric. Then the finite-dimensional distributions of $(\mathbf{Y}_t)$ are regularly varying with index $\alpha$ and the limit $\gamma$ in* (2.4) *exists. If $\gamma > 0$, then*

$$(3.2) \qquad a_n^{-1} \mathbf{S}_n \xrightarrow{d} \mathbf{X}_\alpha,$$

*where the sequence $(a_n)$ is given by*

$$P(|Y_1| > a_n) \sim n^{-1}$$

*and $\mathbf{X}_\alpha$ is an $\alpha$-stable random vector.*

REMARK 3.3. In the case when $E|\mathbf{G}_1|^{2+\delta} + E|\mathbf{Y}_1|^{2+\delta} < \infty$ and $EY_1 = 0$, (3.2) turns into $n^{-1/2}\mathbf{S}_n \overset{d}{\to} \mathbf{X}$, where $\mathbf{X}$ is Gaussian with mean zero and the same covariance structure as $\mathbf{G}_1$. This follows since $(\mathbf{G}_t Y_t)$ is a strictly stationary martingale sequence; see [4].

REMARK 3.4. It is not difficult to see that $\mathbf{Y}_t$ is regularly varying with index $\alpha$. For the proof we need a result of Breiman [11]. It says that if one has two independent random variables $\xi, \eta > 0$ a.s., $\xi$ is regularly varying with index $\alpha > 0$ and $E\eta^\nu < \infty$ for some $\nu > \alpha$, then

$$P(\xi\eta > x) \sim E\eta^\alpha P(\xi > x),$$

that is, $\xi\eta$ is regularly varying with the same index $\alpha$. Now observe that, for $t, x > 0$ and a Borel set $S \subset \mathbb{S}^{d-1}$, by multiple application of Breiman's result,

$$\frac{P(|\mathbf{G}_1||Y_1| > tx, \mathbf{G}_1 Y_1/|\mathbf{G}_1||Y_1| \in S)}{P(|\mathbf{G}_1||Y_1| > x)}$$

$$= \frac{P(|\mathbf{G}_1||Y_1| > tx, \mathrm{sign}(Y_1)\widetilde{\mathbf{G}}_1 \in S)}{P(|\mathbf{G}_1||Y_1| > x)}$$

$$= \frac{P(|\mathbf{G}_1|Y_1 > tx, \widetilde{\mathbf{G}}_1 \in S)}{P(|\mathbf{G}_1||Y_1| > x)} + \frac{P(|\mathbf{G}_1|Y_1 < -tx, -\widetilde{\mathbf{G}}_1 \in S)}{P(|\mathbf{G}_1||Y_1| > x)}$$

$$\sim \frac{E(|\mathbf{G}_1|^\alpha I_S(\widetilde{\mathbf{G}}_1))P(Y_1 > tx)}{E|\mathbf{G}_1|^\alpha P(|Y_1| > x)} + \frac{E(|\mathbf{G}_1|^\alpha I_S(-\widetilde{\mathbf{G}}_1))P(Y_1 \le -tx)}{E|\mathbf{G}_1|^\alpha P(|Y_1| > x)}.$$

Writing for some $p, q \ge 0$ with $p + q = 1$ and a slowly varying function $L(x)$,

$$P(Y_1 > x) = pL(x)x^{-\alpha} \quad \text{and} \quad P(Y_1 \le -x) = qL(x)|x|^{-\alpha}, \qquad x > 0,$$

we can read off the spectral measure of the vector $\mathbf{Y}_1$:

$$(3.3) \qquad P(\mathbf{\Theta} \in S) = p\frac{E(|\mathbf{G}_1|^\alpha I_S(\widetilde{\mathbf{G}}_1))}{E|\mathbf{G}_1|^\alpha} + q\frac{E(|\mathbf{G}_1|^\alpha I_S(-\widetilde{\mathbf{G}}_1))}{E|\mathbf{G}_1|^\alpha}.$$

By regular variation, $a_n = n^{1/\alpha}\ell(n)$ for some slowly varying function $\ell$. By Breiman's result and since $E|\mathbf{G}_1|^{\alpha+\epsilon} < \infty$ for some $\epsilon > 0$, it also follows that

$$P(|\mathbf{G}_1||Y_1| > x) \sim E|\mathbf{G}_1|^\alpha P(|Y_1| > x),$$

and, therefore, $P(|\mathbf{Y}_1| > ca_n) \sim n^{-1}$ for some constant $c > 0$. Moreover, we have

$$(3.4) \qquad nP(a_n^{-1}\mathbf{Y}_1 \in \cdot) \overset{v}{\to} \mu_1,$$

for some measure $\mu_1$ on $\overline{\mathbb{R}}^d \setminus \{\mathbf{0}\}$ which is determined by $\alpha$ and the spectral measure.

REMARK 3.5. It follows from the proof below that

$$nP(a_n^{-1}(\mathbf{Y}_1, \ldots, \mathbf{Y}_h) \in d(\mathbf{x}_1, \ldots, \mathbf{x}_h))$$

(3.5)
$$\overset{v}{\to} \mu_1(d\mathbf{x}_1)\varepsilon_{\mathbf{0}}(d(\mathbf{x}_2, \ldots, \mathbf{x}_h)) + \cdots + \mu_1(d\mathbf{x}_h)\varepsilon_{\mathbf{0}}(d(\mathbf{x}_1, \ldots, \mathbf{x}_{h-1}))$$

$$=: \mu_h(d(\mathbf{x}_1, \ldots, \mathbf{x}_h)),$$

where $\mu_1$ is defined by (3.4), $\varepsilon_{\mathbf{0}}$ is the Dirac measure at $\mathbf{0}$ and

(3.6)
$$(\mathbf{Y}_1, \ldots, \mathbf{Y}_h) := \text{vec}(\mathbf{Y}_1, \ldots, \mathbf{Y}_h) \quad \text{and}$$

$$(\mathbf{x}_1, \ldots, \mathbf{x}_h) := \text{vec}(\mathbf{x}_1, \ldots, \mathbf{x}_h).$$

This means, in particular, that the limiting measure in the definition of regular variation for $(\mathbf{Y}_1, \ldots, \mathbf{Y}_h)$ is the same as in the definition of regular variation for $\text{vec}(\mathbf{Y}'_1, \ldots, \mathbf{Y}'_h)$, where $\mathbf{Y}'_i$ are i.i.d. copies of $\mathbf{Y}_1$. This part of the theorem is valid for any $\alpha > 0$.

PROOF OF THEOREM 2.2. We verify the conditions of Theorem 2.2. Since A.3 implies $\mathcal{A}(a_n)$ and since we require $\gamma > 0$, it remains to check (a) and (c) in Theorem 2.2.

(a) *Regular variation of the finite-dimensional distributions.* We show regular variation of the vector $(\mathbf{Y}_1, \ldots, \mathbf{Y}_h)$ defined in (3.6), that is, we show that (3.5) holds.

We restrict ourselves to proof of regular variation of the pairs $(\mathbf{Y}_1, \mathbf{Y}_2) := \text{vec}(\mathbf{Y}_1, \mathbf{Y}_2)$; the case of general finite-dimensional distributions is completely analogous. The regular variation of $\mathbf{Y}_1$ was explained in Remark 3.4. Let now $B_1$ and $B_2$ be two Borel sets in $[0, \infty]^d \setminus \{\mathbf{0}\}$, bounded away from zero. In particular, there exists $M > 0$ such that $|\mathbf{x}| > M$ for all $\mathbf{x} \in B_1$ and $\mathbf{x} \in B_2$. Then for any $\epsilon > 0$, by intersecting with the events $\{|\mathbf{G}_i| \leq \epsilon\}$ and $\{|\mathbf{G}_i| > \epsilon\}$, $i = 1, 2$,

$$\{a_n^{-1}\mathbf{Y}_1 \in B_1, a_n^{-1}\mathbf{Y}_2 \in B_2\}$$

$$\subset \{|\mathbf{G}_1||Y_1| > Ma_n, |\mathbf{G}_2||Y_2| > Ma_n\}$$

$$\subset \{\epsilon|Y_1| > Ma_n, \epsilon|Y_2| > Ma_n\}$$

$$\cup \{|\mathbf{G}_1|I_{(\epsilon,\infty)}(|\mathbf{G}_1|)|Y_1| > Ma_n, \epsilon|Y_2| > Ma_n\}$$

$$\cup \{|\mathbf{G}_2|I_{(\epsilon,\infty)}(|\mathbf{G}_2|)|Y_2| > Ma_n, \epsilon|Y_1| > Ma_n\}$$

$$\cup \{|\mathbf{G}_1|I_{(\epsilon,\infty)}(|\mathbf{G}_1|)|Y_1| > Ma_n, |\mathbf{G}_2|I_{(\epsilon,\infty)}(|\mathbf{G}_2|)|Y_2| > Ma_n\}$$

$$=: \bigcup_{i=1}^{4} D_i.$$

By independence and an application of Breiman's result, $nP(D_1) \to 0$ and $nP(D_2) \to 0$. Similarly,

$$nP(D_3) \leq nP(|\mathbf{G}_2|I_{(\epsilon,\infty)}(|\mathbf{G}_2|)|Y_2| > Ma_n)$$
$$\sim nP(|Y_2| > Ma_n)E(|\mathbf{G}_2|^\alpha I_{(\epsilon,\infty)}(|\mathbf{G}_2|)),$$

and thus, by Lebesgue's dominated convergence theorem,

$$\lim_{\epsilon \uparrow \infty} \limsup_{n \to \infty} nP(D_3) = 0,$$

and $nP(D_4) \to 0$ can be proved in the same way. We conclude that

$$nP(a_n^{-1}(\mathbf{Y}_1, \mathbf{Y}_2) \in d(\mathbf{x}_1, \mathbf{x}_2)) \xrightarrow{v} \mu_1(d\mathbf{x}_1)\varepsilon_\mathbf{0}(d\mathbf{x}_1) + \mu_1(d\mathbf{x}_2)\varepsilon_\mathbf{0}(d\mathbf{x}_2)$$
$$= \mu_2(d(\mathbf{x}_1, \mathbf{x}_2));$$

see [23]. This proves the regular variation of the two-dimensional finite-dimensional distributions. The higher-dimensional case is completely analogous.

(c) *The condition* (2.9). We have for any $y > 0$,

$$P\left(\max_{k \leq t \leq r_n} |\mathbf{G}_t||Y_t| > ya_n \,\Big|\, |\mathbf{G}_0||Y_0| > ya_n\right)$$
$$\leq P\left(\max_{k \leq t \leq r_n} |\mathbf{G}_t| > ya_n/(s_k a_{r_n}) \,\Big|\, |\mathbf{G}_0||Y_0| > ya_n\right)$$
$$+ P\left(\max_{k \leq t \leq r_n} |Y_t| > s_k a_{r_n}\right)$$
$$=: I_1 + I_2,$$

where $(s_k)$ is any sequence such that $s_k \to \infty$. In what follows all calculations go through for any $y > 0$; for ease of notation, we set $y = 1$. Then, by Remark 3.1,

$$\lim_{k \to \infty} \lim_{n \to \infty} I_2 = \lim_{k \to \infty}(1 - \Phi_\alpha(s_k)) = 0.$$

An application of Markov's inequality yields, for some constant $c > 0$ and $\epsilon > 0$ as in A.2 (here and in what follows, $c$ denotes any positive constant whose value is not of interest),

$$I_1 \leq \sum_{t=k}^{r_n} P(|\mathbf{G}_t| > a_n/(s_k a_{r_n}) \,|\, |\mathbf{G}_0||Y_0| > a_n)$$
$$\leq \left(\frac{s_k a_{r_n}}{a_n}\right)^{\alpha+\epsilon} \sum_{t=k}^{r_n} \frac{E[|\mathbf{G}_t|^{\alpha+\epsilon}I_{\{|\mathbf{G}_0||Y_0| > a_n/(s_k a_{r_n})\}}]}{P(|\mathbf{G}_0||Y_0| > a_n)}$$

$$\leq cnr_n\left(\frac{s_k a_{r_n}}{a_n}\right)^{\alpha+\epsilon} E|\mathbf{G}_0|^{\alpha+\epsilon}$$

$$\to 0 \qquad \text{as } n \to \infty.$$

Here we used Breiman's result [11] to show that

$$P(|\mathbf{G}_0||Y_0| > a_n) \sim E|\mathbf{G}_0|^\alpha P(|Y_0| > a_n),$$

condition (3.1) and the fact that $E|\mathbf{G}_1|^{\alpha+\epsilon} < \infty$; see A.2.

Now we turn to

$$P\bigg(\max_{-r_n \leq t \leq -k} |\mathbf{G}_t||Y_t| > a_n \,\Big|\, |\mathbf{G}_0||Y_0| > a_n\bigg)$$

$$\leq P\bigg(\max_{-r_n \leq t \leq -k} |\mathbf{G}_t| > a_n/(s_k a_{r_n}) \,\Big|\, |\mathbf{G}_0||Y_0| > a_n\bigg)$$

$$+ P\bigg(\max_{-r_n \leq t \leq -k} |Y_t| > s_k a_{r_n} \,\Big|\, |\mathbf{G}_0||Y_0| > a_n\bigg)$$

$$=: I_3 + I_4.$$

The quantity $I_3$ can be treated in the same way as $I_1$ to show that $I_3 \to 0$ a.s. as $n \to \infty$. We turn to $I_4$. Fix $0 < M < \infty$. Then

$$I_4 \leq \frac{P(\max_{-r_n \leq t \leq -k} |Y_t| > s_k a_{r_n}, M|Y_0| > a_n)}{P(|\mathbf{G}_0||Y_0| > a_n)}$$

$$+ \frac{P(|\mathbf{G}_0| I_{(M,\infty)}(|\mathbf{G}_0|)|Y_0| > a_n)}{P(|\mathbf{G}_0||Y_0| > a_n)}$$

$$=: I_{41} + I_{42}.$$

By independence of the $Y_i$'s, Breiman's [11] result and since $r_n \to \infty$,

$$I_{41} \sim \frac{P(\max_{-r_n \leq t \leq -k} |Y_t| > s_k a_{r_n}) M^\alpha P(|Y_0| > a_n)}{E|\mathbf{G}_0|^\alpha P(|Y_0| > a_n)}$$

$$\sim c(1 - \Phi_\alpha(s_k)) \qquad \text{as } n \to \infty$$

$$\to 0 \qquad \text{as } k \to \infty.$$

By virtue of Breiman's [11] result,

$$I_{42} \sim \frac{E(|\mathbf{G}_0|^\alpha I_{(M,\infty)}(|\mathbf{G}_0|))P(|Y_0| > a_n)}{E|\mathbf{G}_0|^\alpha P(|Y_0| > a_n)}.$$

Since $|\mathbf{G}_0|$ has finite moments of order greater than $\alpha$, an application of the Lebesgue dominated convergence theorem yields

$$\lim_{M \to \infty} \lim_{n \to \infty} I_{42} = 0.$$

This proves (2.9).   □

Thus, the conditions (a)–(c) and $\gamma > 0$ of Theorem 2.2 are satisfied. In the case $\alpha < 1$, Theorem 2.2 immediately yields (3.2). In the case $\alpha \in [1, 2)$, we have to check condition (2.11). It suffices to show it for components $\mathbf{S}_n^{(i)}(0, y]$, $i = 1, \ldots, d$, of $\mathbf{S}_n(0, y]$. Since the components can be handled in the same way, we suppress the dependence on $i$ and, for ease of notation, write $G_t Y_t$ for the summands of the $i$th component.

We start with the case $\alpha \in (1, 2)$. As before, write $\mathcal{F}_t = \sigma(Y_{t-1}, Y_{t-2}, \ldots)$. Then, for $z > 0$, since $EY_1 = 0$,

$$E[G_t Y_t I_{(0,z]}(|G_t Y_t|/a_n) \mid \mathcal{F}_t] = G_t E[Y_t I_{(0,z]}(|G_t Y_t|/a_n) \mid G_t]$$
$$= -G_t E[Y_t I_{(z,\infty)}(|G_t Y_t|/a_n) \mid G_t].$$

Consider the decomposition

$$a_n^{-1} \sum_{t=1}^n [G_t Y_t I_{(0,z]}(|G_t Y_t|/a_n) - E[G_1 Y_1 I_{(0,z]}(|G_1 Y_1|/a_n)]]$$

$$= a_n^{-1} \sum_{t=1}^n [G_t Y_t I_{(0,z]}(|G_t Y_t|/a_n) - G_t E[Y_t I_{(0,z]}(|G_t Y_t|/a_n) \mid G_t]]$$

$$- a_n^{-1} \sum_{t=1}^n [G_t E[Y_t I_{(z,\infty)}(|G_t Y_t|/a_n) \mid G_t] - E[G_1 Y_1 I_{(z,\infty)}(|G_1 Y_1|/a_n)]]$$

$$=: T_1 + T_2.$$

For fixed $n$, $T_1$ is a sum of stationary mean zero martingale differences. An application of Karamata's theorem ([5], page 26) to the regularly varying random variable $G_1 Y_1$ with index $\alpha$ yields for some constant $c > 0$,

$$
\begin{aligned}
\operatorname{var}(T_1) &= n a_n^{-2} E[G_1 Y_1 I_{(0,z]}(|G_1 Y_1|/a_n) \\
&\qquad - G_1 E[Y_1 I_{(0,z]}(|G_1 Y_1|/a_n) \mid G_1]]^2 \\
&\leq c n a_n^{-2} E[G_1 Y_1 I_{(0,z]}(|G_1 Y_1|/a_n)]^2 \\
&\sim c z^{2-\alpha} \qquad \text{as } n \to \infty \\
&\to 0 \qquad \text{as } z \downarrow 0.
\end{aligned}
$$
(3.7)

Next we treat $T_2$. Fix $0 < \delta < M < \infty$ to be chosen later. Notice that, by Karamata's theorem and the uniform convergence theorem for regularly varying functions uniformly for $c \in [\delta, M]$,

$$\frac{E[Y_1 I_{(cx,\infty)}(|Y_1|)]}{cx P(|Y_1| > cx)} \to C$$

for some constant $C$. Taking this into account, the strong law of large numbers yields, with probability 1,

$$
\begin{aligned}
a_n^{-1} &\sum_{t=1}^{n} G_t I_{[\delta,M]}(|G_t|) E[Y_t I_{(z,\infty)}(|G_t Y_t|/a_n) \mid G_t] \\
&= a_n^{-1} \sum_{t=1}^{n} G_t I_{[\delta,M]}(|G_t|) \\
&\qquad \times [(za_n/G_t) P(|Y_t| > za_n/|G_t| \mid G_t)(C + o(1))] \\
&= (C + o(1)) z^{1-\alpha} n^{-1} \sum_{t=1}^{n} |G_t|^\alpha I_{[\delta,M]}(|G_t|) \\
&\to C z^{1-\alpha} E[|G_1|^\alpha I_{[\delta,M]}(|G_1|)].
\end{aligned}
$$
(3.8)

On the other hand, since $G_1 I_{[\delta,M]}(|G_1|) Y_1$ is regularly varying with index $\alpha \in (1,2)$, by the same argument and Breiman's result,

$$
\begin{aligned}
n a_n^{-1} &E[G_1 I_{[\delta,M]}(|G_1|) Y_1 I_{(z,\infty)}(|G_1 Y_1|/a_n)] \\
&= n a_n^{-1}[(C + o(1))(za_n) P(G_1 I_{[\delta,M]}(|G_1|)|Y_1| > za_n)] \\
&= (C + o(1)) z^{1-\alpha} E[|G_1|^\alpha I_{[\delta,M]}(|G_1|)].
\end{aligned}
$$
(3.9)

This shows that (3.8) and (3.9) cancel asymptotically as $n \to \infty$ for every fixed $z$.

A similar argument shows that, with probability 1,

$$
\begin{aligned}
a_n^{-1} &\left| \sum_{t=1}^{n} G_t I_{[0,\delta]}(|G_t|) E[Y_t I_{(z,\infty)}(|G_t Y_t|/a_n) \mid G_t] \right| \\
&\leq a_n^{-1} \sum_{t=1}^{n} |G_t| I_{[0,\delta]}(|G_t|) E[|Y_1| I_{(z,\infty)}(\delta|Y_1|/a_n)] \\
&\to c(z/\delta)^{1-\alpha} E[|G_1| I_{[0,\delta]}(|G_1|)].
\end{aligned}
$$
(3.10)

Moreover,

$$
\begin{aligned}
n a_n^{-1} &|E[G_1 I_{[0,\delta]}(|G_1|) Y_1 I_{(z,\infty)}(|G_1 Y_1|/a_n)]| \\
&\leq n a_n^{-1} E[|G_1| I_{[0,\delta]}(|G_1|)|Y_1| I_{(z,\infty)}(\delta|Y_1|/a_n)] \\
&\sim c(z/\delta)^{1-\alpha} E[|G_1| I_{[0,\delta]}(|G_1|)].
\end{aligned}
$$
(3.11)

Now choose $\delta = z^2$. Then, first letting $n \to \infty$ and then $z \downarrow 0$, both (3.10) and (3.11) vanish asymptotically.

Finally, we consider

$$a_n^{-1} E \left| \sum_{t=1}^n G_t I_{(M,\infty)}(|G_t|) E[Y_t I_{(z,\infty)}(|G_t Y_t|/a_n) \mid G_t] \right|$$

$$\leq a_n^{-1} n E[|G_1| I_{(M,\infty)}(|G_1|) |Y_1| I_{(z,\infty)}(|G_1 Y_1|/a_n)].$$

An application of Breiman's result to the regularly varying random variable $G_1 I_{[M,\infty)}(|G_1|) Y_1$ gives that the right-hand side is asymptotically equivalent as $n \to \infty$ to

$$c z^{1-\alpha} E[|G_1|^\alpha I_{[M,\infty)}(|G_1|)].$$

Choosing $M$ large enough, the right-hand side is smaller than $z$, say. The same argument can be applied to

$$n a_n^{-1} |E[G_1 I_{[M,\infty)}(|G_1|) Y_1 I_{(z,\infty)}(|G_1 Y_1|/a_n)]|.$$

Collecting the bounds above, we see that

$$\lim_{z \downarrow 0} \limsup_{n \to \infty} P(|T_2| > r) = 0, \qquad r > 0.$$

This together with (3.7) concludes the proof of (2.11) for $\alpha \in (1, 2)$.

For $\alpha = 1$, we use the additional condition of symmetry of $Y_t$. Then $E\mathbf{S}_n(0, y] = 0$ and the same argument as for $\text{var}(T_1)$ above shows that (2.11) holds in this case as well. This concludes the proof of (2.11).

Since the conditions of Theorem 2.2 are satisfied for $\alpha \in [1, 2)$, we conclude that

$$a_n^{-1}(\mathbf{S}_n - E\mathbf{S}_n(0, 1]) \xrightarrow{d} \mathbf{X}_\alpha$$

for some $\alpha$-stable random vector in $\mathbb{R}^d$. For $\alpha = 1$, we can drop $E\mathbf{S}_n(0, y]$ because of the symmetry of $\mathbf{G}_t Y_t$. For $\alpha \in (1, 2)$, $\mathbf{G}_t Y_t$ is regularly varying with index $\alpha$. Since $E(\mathbf{G}_t Y_t) = 0$, Karamata's theorem yields

$$a_n^{-1} E\mathbf{S}_n(0, 1] \to \mathbf{b}$$

for some constant $\mathbf{b}$ which can be incorporated in the stable limit, and, therefore, centering in (3.2) can be avoided. This concludes the proof of Theorem 3.2.

**4. Gaussian quasi maximum likelihood estimation for GARCH processes with heavy-tailed innovations.** In this section we apply Theorem 3.2 to Gaussian quasi maximum likelihood estimation (QMLE) in GARCH processes. The limit properties of the QMLE were studied by Berkes et al. [3]. They proved strong consistency of the QMLE under the moment condition $E|Z_1|^{2+\delta} < \infty$ for some $\delta > 0$ and established asymptotic normality under $EZ_1^4 < \infty$. Here $(Z_t)$ is an i.i.d. innovation sequence; see Section 4.1 below for the definition of the GARCH model and the QMLE. Hall and Yao [16] refined these results and also allowed for innovations sequences, where $Z_1^2$ is regularly varying with index $\alpha \in (1,2)$. Then the speed of convergence is slower than the usual $\sqrt{n}$ rate and the limiting distribution of the QMLE is (multivariate) $\alpha$-stable.

It is our objective to show that the asymptotic theories for the QMLE under light- and heavy-tailed innovations parallel each other and that very similar techniques can be applied in both cases. However, in the light-tailed case (see [3]) an application of the CLT for stationary ergodic martingale differences is the basic tool which establishes the asymptotic normality of the QMLE. In the heavy-tailed situation one depends on an analog of the CLT which is provided by Theorem 3.2.

As a matter of fact, the structure of the proofs shows that the asymptotic properties of the QMLE are not dependent on the particular structure of the GARCH process if one can establish the regular variation of the finite-dimensional distributions of the underlying process $(X_t)$ and the mixing condition $\mathcal{A}(a_n)$. Therefore, the results of this section have the potential to be extended to more general models, including, for example, the AGARCH or EGARCH models whose QMLE properties in the light-tailed case are treated in [28]. The most intricate step in the proof is, however, the verification of this mixing condition for a given time series model. We establish this condition for a GARCH process by an adaptation of Theorem 4.3 in [21]; this yields strong mixing with geometric rate of the relevant sequence. We devote Section 4.4 to the solution of this problem.

Before we start, we introduce some notation. If $K \subset \mathbb{R}^d$ is a compact set, we write $\mathbb{C}(K, \mathbb{R}^{d'})$ for the space of continuous $\mathbb{R}^{d'}$-valued functions equipped with the sup-norm $\|v\|_K = \sup_{s \in K} |v(s)|$. The space $\mathbb{C}(K, \mathbb{R}^{d_1 \times d_2})$ consists of the continuous $d_1 \times d_2$-matrix valued functions on $K$; in $\mathbb{R}^{d_1 \times d_2}$ we work with the operator norm induced by the Euclidean norm $|\cdot|$, that is,

$$\|\mathbf{A}\| = \sup_{|x|=1} |\mathbf{A}x|, \qquad \mathbf{A} \in \mathbb{R}^{d_1 \times d_2}.$$

### 4.1. *Definition of the QMLE.*

Recall the definition of a GARCH$(p,q)$ process $(X_t)$ from (1.1). As before, $(Z_t)$ is an i.i.d. innovation sequence with $EZ_1^2 = 1$ and $EZ_1 = 0$, and $\alpha_i, \beta_j$ are nonnegative constants. GARCH processes have been intensively investigated over the last few years. Assumptions for strict stationarity are complicated: they are expressed in terms of Lyapunov exponents of certain random matrices; see [6] for details. A necessary condition for stationarity is

$$(4.1) \qquad \beta_1 + \cdots + \beta_q < 1$$

(Corollary 2.3 in [6]). We will make use of this condition later.

In what follows we always assume strict stationarity of the GARCH processes. As a matter of fact, the observation $X_t$ is always a measurable function of the past and present innovations $(Z_t, Z_{t-1}, Z_{t-2}, \ldots)$; hence, $(X_t)$ is automatically ergodic.

In what follows we review how an approximation to the conditional *Gaussian* likelihood of a stationary GARCH$(p, q)$ process is constructed, that is, a conditional likelihood under the *synthetic* assumption $Z_t$ i.i.d. $\sim \mathcal{N}(0, 1)$. Given $X_0, \ldots, X_{-p+1}$ and $\sigma_0^2, \ldots, \sigma_{-q+1}^2$, the random variables $X_1, \ldots, X_n$ are conditionally Gaussian with mean zero and variances $h_t(\boldsymbol{\theta})$, $t = 1, \ldots, n$, where $\boldsymbol{\theta} = (\alpha_0, \alpha_1, \ldots, \alpha_p, \beta_1, \ldots, \beta_q)^T$ denotes the presumed parameter and

$$\check{h}_t(\boldsymbol{\theta}) = \begin{cases} \sigma_t^2, & t \leq 0, \\ \alpha_0 + \alpha_1 X_{t-1}^2 + \cdots + \alpha_p X_{t-p}^2 \\ \quad + \beta_1 \check{h}_{t-1}(\boldsymbol{\theta}) + \cdots + \beta_q \check{h}_{t-q}(\boldsymbol{\theta}), & t > 0. \end{cases}$$

The conditional Gaussian log-likelihood has the form

(4.2)
$$\log f_{\boldsymbol{\theta}}(X_1, \ldots, X_n \mid X_0, \ldots, X_{-p+1}, \sigma_0^2, \ldots, \sigma_{-q+1}^2)$$
$$= -\frac{n}{2} \log(2\pi) - \frac{1}{2} \sum_{t=1}^{n} \left( \frac{X_t^2}{\check{h}_t(\boldsymbol{\theta})} + \log \check{h}_t(\boldsymbol{\theta}) \right).$$

Since $X_0, \ldots, X_{-p+1}$ are not available and the squared volatilities $\sigma_0^2, \ldots, \sigma_{-q+1}^2$ unobservable, the conditional Gaussian log-likelihood (4.2) cannot be numerically evaluated without a certain initialization for $\sigma_0^2, \ldots, \sigma_{-p+1}^2$ and $X_0, \ldots, X_{-q+1}$. The initial values being asymptotically irrelevant, we set the $X_t$'s equal to zero and $\hat{h}_t(\boldsymbol{\theta}) = \alpha_0/(1 - \beta_1 - \cdots - \beta_q)$ for $t \leq 0$. We arrive at

(4.3)  $\hat{h}_t(\boldsymbol{\theta}) = \begin{cases} \alpha_0/(1 - \beta_1 - \cdots - \beta_q), & t \leq 0, \\ \alpha_0 + \alpha_1 X_{t-1}^2 + \cdots + \alpha_{\min(p, t-1)} X_{\max(t-p, 1)}^2 \\ \quad + \beta_1 \hat{h}_{t-1}(\boldsymbol{\theta}) + \cdots + \beta_q \hat{h}_{t-q}(\boldsymbol{\theta}), & t > 0. \end{cases}$

The function $(\hat{h}_t(\boldsymbol{\theta}))^{1/2}$ can be understood as an estimate of the volatility at time $t$ and under parameter hypothesis $\boldsymbol{\theta}$. It can be established that $|\hat{h}_t - \check{h}_t| \overset{\text{a.s.}}{\to} 0$ with a geometric rate of convergence and uniformly on the compact set $K$ defined in (4.4) below. This suggests that, by replacing $\check{h}_t(\boldsymbol{\theta})$ by $\hat{h}_t(\boldsymbol{\theta})$ in (4.2), we obtain a good approximation to the conditional Gaussian log-likelihood. Since the constant $-n \log(2\pi)/2$ does not matter for the optimization, we define the QMLE $\hat{\boldsymbol{\theta}}_n$ as a maximizer of the function

$$\widehat{L}_n(\boldsymbol{\theta}) = \sum_{t=1}^{n} \hat{\ell}_t(\boldsymbol{\theta}) = -\frac{1}{2} \sum_{t=1}^{n} \left( \frac{X_t^2}{\hat{h}_t(\boldsymbol{\theta})} + \log \hat{h}_t(\boldsymbol{\theta}) \right)$$

with respect to $\boldsymbol{\theta} \in K$, with $K$ being the compact set

(4.4)      $K = \{ \boldsymbol{\theta} \in \mathbb{R}^{p+q+1} \mid m \leq \alpha_i, \beta_j \leq M, \beta_1 + \cdots + \beta_q \leq \bar{\beta} \},$

where $0 < m < M < \infty$ and $0 < \bar{\beta} < 1$ are such that $qm < \bar{\beta}$.

REMARK 4.1. From a comparison with [3], one might think at first sight that our definition of the QMLE is different from theirs. To see that $\hat{h}_t$ coincides with $\tilde{w}_t$ in [3], introduce the polynomials

$$\boldsymbol{\alpha}(z) = \alpha_1 z + \cdots + \alpha_p z^p \quad \text{and} \quad \boldsymbol{\beta}(z) = 1 - \beta_1 z - \cdots - \beta_q z^q$$

for every $\boldsymbol{\theta} = (\alpha_0, \alpha_1, \ldots, \alpha_p, \beta_1, \ldots, \beta_q)^T \in K$. Then one can show by induction on $t$ that

$$(4.5) \qquad \hat{h}_t(\boldsymbol{\theta}) = \frac{\alpha_0}{\boldsymbol{\beta}(1)} + \sum_{j=1}^{t-1} \psi_j(\boldsymbol{\theta}) X_{t-j}^2,$$

where the coefficients $\psi_j(\boldsymbol{\theta})$ are defined through

$$(4.6) \qquad \frac{\boldsymbol{\alpha}(z)}{\boldsymbol{\beta}(z)} = \sum_{j=1}^{\infty} \psi_j(\boldsymbol{\theta}) z^j, \qquad |z| \le 1.$$

Note that the latter Taylor series representation is valid because $\beta_i \ge 0$ and $\beta_1 + \cdots + \beta_q \le \bar{\beta} < 1$ imply $\boldsymbol{\beta}(z) \ne 0$ on $K$ for $|z| \le 1 + \epsilon$ and $\epsilon > 0$ sufficiently small. We choose (4.3) rather than (4.5) as a first definition for the squared volatility estimate under parameter hypothesis $\boldsymbol{\theta}$, because the recursion (4.3) is natural and computationally attractive. In [3] the starting point for the definition of the QMLE is Theorem 2.2, which says that for all $t \in \mathbb{Z}$ one has $h_t(\boldsymbol{\theta}_0) = \sigma_t^2$, where $\boldsymbol{\theta}_0$ is the true parameter and

$$(4.7) \qquad h_t(\boldsymbol{\theta}) = \frac{\alpha_0}{\boldsymbol{\beta}(1)} + \sum_{j=1}^{\infty} \psi_j(\boldsymbol{\theta}) X_{t-j}^2.$$

In [3] this leads to the definition of a squared volatility estimate at time $t$ under parameter $\boldsymbol{\theta}$ based on $(X_1, \ldots, X_n)$, which is given by (4.5). Note also that $(h_t(\boldsymbol{\theta}))$ obeys

$$(4.8) \qquad \begin{aligned} h_{t+1}(\boldsymbol{\theta}) &= \alpha_0 + \alpha_1 X_t^2 + \cdots + \alpha_p X_{t+1-p}^2 \\ &\quad + \beta_1 h_t(\boldsymbol{\theta}) + \cdots + \beta_q h_{t+1-q}(\boldsymbol{\theta}), \qquad \boldsymbol{\theta} \in K. \end{aligned}$$

4.2. *Limit distribution in the case* $EZ_1^4 < \infty$. First we list the conditions employed by [3] for establishing consistency and asymptotic normality of $\hat{\boldsymbol{\theta}}_n$. Write $\boldsymbol{\theta}_0 = (\alpha_0^\circ, \alpha_1^\circ, \ldots, \alpha_p^\circ, \beta_1^\circ, \ldots, \beta_q^\circ)^T$ for the true parameter.

C.1. There is a $\delta > 0$ such that $E|Z_1|^{2+\delta} < \infty$.

C.2. The distribution of $|Z_1|$ is not concentrated in one point.

C.3. There is a $\mu > 0$ such that $P(|Z_1| \le t) = o(t^\mu)$ as $t \downarrow 0$.

C.4. The true parameter $\boldsymbol{\theta}_0$ lies in the interior of $K$.

C.5. The polynomials $\boldsymbol{\alpha}^\circ(z) = \alpha_1^\circ z + \cdots + \alpha_p^\circ z^p$ and $\boldsymbol{\beta}^\circ(z) = 1 - \beta_1^\circ z - \cdots - \beta_q^\circ z^q$ do not have any common roots.

Now we are ready to quote the main result of [3]. We cite it in order to be able to compare the assumptions and assertions both in the light- and heavy-tailed cases; see Theorem 4.4 below.

THEOREM 4.2 (Theorem 4.1 of [3]).   *Let* $(X_t)$ *be a stationary* GARCH$(p, q)$ *process with true parameter vector* $\boldsymbol{\theta}_0$. *Suppose the conditions* C.1–C.5 *hold. Then the QMLE* $\hat{\boldsymbol{\theta}}_n$ *is strongly consistent, that is,*

$$\hat{\boldsymbol{\theta}}_n \overset{\text{a.s.}}{\to} \boldsymbol{\theta}_0, \qquad n \to \infty.$$

*If, in addition,* $EZ_0^4 < \infty$, *then* $\hat{\boldsymbol{\theta}}_n$ *is also asymptotically normal, that is,*

$$\sqrt{n}(\hat{\boldsymbol{\theta}}_n - \boldsymbol{\theta}_0) \overset{d}{\to} \mathcal{N}(0, \mathbf{B}_0^{-1} \mathbf{A}_0 \mathbf{B}_0^{-1}),$$

*where the* $(p + q + 1) \times (p + q + 1)$ *matrices* $\mathbf{A}_0$ *and* $\mathbf{B}_0$ *are given by*

$$
\begin{aligned}
(4.9) \qquad \mathbf{A}_0 &= \frac{E(Z_0^4 - 1)}{4} E\left(\frac{1}{\sigma_1^4} h_1'(\boldsymbol{\theta}_0)^T h_1'(\boldsymbol{\theta}_0)\right), \\
\mathbf{B}_0 &= -\frac{1}{2} E\left(\frac{1}{\sigma_1^4} h_1'(\boldsymbol{\theta}_0)^T h_1'(\boldsymbol{\theta}_0)\right).
\end{aligned}
$$

4.3. *Limit distribution in the case* $EZ_1^4 = \infty$.   First we identify the limit determining term for the QMLE. To this end, we set analogously to [3],

$$L_n(\boldsymbol{\theta}) = \sum_{t=1}^{n} \ell_t(\boldsymbol{\theta}) = -\frac{1}{2} \sum_{t=1}^{n} \left(\frac{X_t^2}{h_t(\boldsymbol{\theta})} + \log h_t(\boldsymbol{\theta})\right)$$

and define $\tilde{\boldsymbol{\theta}}_n$ as a maximizer of $L_n$ with respect to $\boldsymbol{\theta} \in K$. It is a slightly simpler problem to analyze $\tilde{\boldsymbol{\theta}}_n$ because $(\ell_t)$ is stationary ergodic, in contrast to $(\hat{\ell}_t)_{t \in \mathbb{N}}$. As is shown in Proposition 4.3 below, $\hat{\boldsymbol{\theta}}_n$ and $\tilde{\boldsymbol{\theta}}_n$ are asymptotically equivalent. It turns out that the asymptotic distribution of the QMLE is essentially determined by the limit behavior of $L_n'(\boldsymbol{\theta}_0)/n$, up to multiplication with the matrix $-\mathbf{B}_0^{-1}$. These results follow by a careful analysis of the proofs in [3]. We omit details and refer to the website [20] for a detailed proof. Compare also with the similar reference [28], where the case of processes with a more general volatility structure than GARCH is treated.

PROPOSITION 4.3.   *Let* $(X_t)$ *be a stationary* GARCH$(p, q)$ *process with true parameter vector* $\boldsymbol{\theta}_0$. *Suppose the conditions* C.1–C.5 *apply. If there is a positive sequence* $(x_n)_{n \geq 1}$ *with* $x_n = o(n)$ *as* $n \to \infty$ *and*

$$(4.10) \qquad x_n \frac{L_n'(\boldsymbol{\theta}_0)}{n} \overset{d}{\to} \mathbf{D}, \qquad n \to \infty,$$

*for an $\mathbb{R}^{p+q+1}$-valued random variable $\mathbf{D}$, then the QMLE $\hat{\boldsymbol{\theta}}_n$ satisfies the limit relation*

$$(4.11) \qquad x_n(\hat{\boldsymbol{\theta}}_n - \boldsymbol{\theta}_0) \xrightarrow{d} -\mathbf{B}_0^{-1}\mathbf{D},$$

*where $\mathbf{B}_0$ is given by* (4.9).

Now we can state the main theorem of this section. We note once again that Hall and Yao [16] derived the identical result by means of different techniques.

THEOREM 4.4. *Let $(X_t)$ be a stationary* GARCH$(p,q)$ *process with true parameter vector $\boldsymbol{\theta}_0$. Suppose that $Z_1^2$ is regularly varying with index $\alpha \in (1,2)$ and that C.3–C.5 hold. Moreover, assume that $Z_1$ has a Lebesgue density $f$, where the closure of the interior of the support $\{f > 0\}$ contains the origin. Define $(x_n) = (na_n^{-1})$, where*

$$P(Z_1^2 > a_n) \sim n^{-1}, \qquad n \to \infty.$$

*Then the QMLE $\hat{\boldsymbol{\theta}}_n$ is consistent and*

$$(4.12) \qquad x_n(\hat{\boldsymbol{\theta}}_n - \boldsymbol{\theta}_0) \xrightarrow{d} \mathbf{D}_\alpha, \qquad n \to \infty,$$

*for some nondegenerate $\alpha$-stable vector $\mathbf{D}_\alpha$.*

Before proving the theorem, we discuss its practical consequences for parameter inference:

- The rate of convergence $x_n$ has—roughly speaking—magnitude $n^{1-1/\alpha}$, which is less than $\sqrt{n}$. The heavier the tails of the innovations, that is, the smaller $\alpha$, the slower is the convergence of $\hat{\boldsymbol{\theta}}_n$ toward the true parameter $\boldsymbol{\theta}_0$.
- The limit distribution of the standardized differences $(\hat{\boldsymbol{\theta}}_n - \boldsymbol{\theta}_0)$ is $\alpha$-stable and, hence, non-Gaussian. The exact parameters of this $\alpha$-stable limit are not explicitly known.
- Confidence bands based on the normal approximation of Theorem 4.2 are false if $EZ_1^4 = \infty$.
- By the definition of a GARCH process, the distribution of the innovations $Z_t$ is unknown. Therefore, assumptions about the heaviness of the tails of its distribution are purely hypothetical. As a matter of fact, the tails of the distribution of $X_t$ can be regularly varying even if $Z_t$ has light tails, such as for the normal distribution; see [2]. Depending on the assumptions on the distribution of $Z_1$, one can develop different asymptotic theories for QMLE of GARCH processes: asymptotic normality as provided by Theorem 4.2 or infinite variance stable distributions as provided by Theorem 4.4.

PROOF OF THEOREM 4.4.    The proof follows by combining Theorem 3.2 and Proposition 4.3. Indeed, setting

$$\mathbf{G}_t = h_t'(\boldsymbol{\theta}_0)/\sigma_t^2, \qquad Y_t = (Z_t^2 - 1)/2 \quad \text{and} \quad \mathbf{Y}_t = \mathbf{G}_t Y_t,$$

one recognizes that

$$(4.13) \qquad L_n'(\boldsymbol{\theta}_0) = \frac{1}{2}\sum_{t=1}^\infty \frac{h_t'(\boldsymbol{\theta}_0)}{\sigma_t^2}(Z_t^2 - 1) = \sum_{t=1}^n \mathbf{G}_t Y_t$$

is a martingale transform. Regular variation of $Z_1^2$ with index $\alpha \in (1, 2)$ implies A.1, but also C.1 and C.2. Condition A.2 is fulfilled because $\|h_1'/h_1\|_K$ has finite moments of any order (Lemma 5.2 of [3]), and so has $\|\mathbf{G}_1\|$. The condition A.3 holds if we can show that $(\mathbf{Y}_t)$ is strongly mixing with geometric rate, in which case we choose $r_n = n^\delta$ in $\mathcal{A}(a_n)$ for any small $\delta > 0$, so that (3.1) immediately follows. This choice of $(r_n)$ is justified by the arguments given in [2]. The strong mixing condition with geometric rate of $(\mathbf{Y}_t)$ will be verified in Section 4.4.

Finally, we have to give an argument for $\gamma > 0$. The latter quantity has interpretation as the extremal index of the sequence $(|\mathbf{Y}_t|)$; see Remark 2.4. According to Theorem 3.7.2 in [18], if $\gamma = 0$ and for some sequence $(u_n)$ the relation $\liminf_{n\to\infty} P(\widetilde{M}_n \le u_n) > 0$ holds, then one neccessarily has $\lim_{n\to\infty} P(M_n \le u_n) = 1$. Here $M_n = \max(|\mathbf{Y}_1|, \ldots, |\mathbf{Y}_n|)$ and $(\widetilde{M}_n)$ is the corresponding sequence of partial maxima for an i.i.d. sequence $(R_i)$, where $R_1$ has the same distribution as $|\mathbf{Y}_1|$.

We want to show by contradiction that $\gamma = 0$, using the above result. The random variable $|\mathbf{Y}_1| \overset{d}{=} R_i$ is regularly varying with index $\alpha$ since $\mathbf{Y}_1$ is regularly varying with index $\alpha$. Hence, $(a_n^{-1}\widetilde{M}_n)$ has a Fréchet limit distribution $\Phi_\alpha(x) = \exp\{-x^{-\alpha}\}$, $x > 0$; see Remark 3.1.

On the other hand, we will show that $P(M_n \le xa_n) \to 1$ does not hold for any positive $x$, thus contradicting the hypothesis $\gamma = 0$. Indeed, straightforward arguments exploiting

$$\sum_{j=1}^\infty \frac{\partial \psi_j(\boldsymbol{\theta})}{\partial \alpha_i} z^j = \frac{z^i}{\boldsymbol{\beta}(z)}, \qquad |z| \le 1,$$

for all $i = 1, \ldots, p$, show that

$$(4.14) \qquad \frac{\partial h_t(\boldsymbol{\theta})}{\partial \alpha_i} \ge 0 \qquad \text{for all } i = 0, \ldots, p,$$

and

$$(4.15) \qquad \sum_{i=0}^p \alpha_i \frac{\partial h_t(\boldsymbol{\theta})}{\partial \alpha_i} = h_t(\boldsymbol{\theta}).$$

Since the Euclidean norm is equivalent to the 1-norm $|\mathbf{x}| = \sum_{i=1}^{p+q+1} |x_i|$ and $\alpha_i \leq M$ on $K$, there is a $c > 0$ such that

$$\frac{|h_t'(\boldsymbol{\theta})|}{h_t(\boldsymbol{\theta})} \geq \frac{c}{h_t(\boldsymbol{\theta})} \sum_{i=0}^{p} \alpha_i \left| \frac{\partial h_t(\boldsymbol{\theta})}{\partial \alpha_i} \right| = \frac{c}{h_t(\boldsymbol{\theta})} \sum_{i=0}^{p} \alpha_i \frac{\partial h_t(\boldsymbol{\theta})}{\partial \alpha_i} = c.$$

Note that the last two equalities in the latter display are a consequence of (4.14) and (4.15). In particular, we proved that $|\mathbf{G}_i| \geq c$ for all $i$ and therefore

$$P(M_n \leq x a_n) = P\left( \max_{i=1,\ldots,n} |\mathbf{G}_i| |Y_i| \leq x a_n \right)$$

$$\leq P\left( \max_{i=1,\ldots,n} |Y_i| \leq c^{-1} x a_n \right).$$

The same classical limit result for maxima as above ensures that the right-hand side probability converges to a Fréchet limit and is never equal to 1 for all positive $x$. Thus, we have proved $\gamma > 0$.

Now, all conditions of Theorem 3.2 are verified so that

$$2 a_n^{-1} L_n'(\boldsymbol{\theta}_0) = 2 x_n \frac{L_n'(\boldsymbol{\theta}_0)}{n} \xrightarrow{d} \widetilde{\mathbf{D}}_\alpha,$$

where $\widetilde{\mathbf{D}}_\alpha$ is $\alpha$-stable [notice that $P((Z_0^2 - 1)/2 > a_n/2) \sim P(Z_0^2 > a_n) \sim n^{-1}$]. Since $x_n/n = a_n^{-1} \to 0$, Proposition 4.3 implies

$$x_n(\hat{\boldsymbol{\theta}}_n - \boldsymbol{\theta}_0) \xrightarrow{d} -2^{-1} \mathbf{B}_0^{-1} \widetilde{\mathbf{D}}_\alpha = \mathbf{D}_\alpha.$$

Recalling that a linear transformation of an $\alpha$-stable random vector is again $\alpha$-stable (see [27]), we conclude the proof of the theorem. $\square$

4.4. *Verification of strong mixing with geometric rate of* $(\mathbf{Y}_t)$. To begin with, we quote a powerful result due to Mokkadem [21], which allows one to establish strong mixing in stationary solutions of so-called polynomial linear stochastic recurrence equations (SREs). A sequence $(\mathbf{Y}_t)$ of random vectors in $\mathbb{R}^d$ obeys a linear SRE if

$$(4.16) \qquad \mathbf{Y}_t = \mathbf{P}_t \mathbf{Y}_{t-1} + \mathbf{Q}_t,$$

where $((\mathbf{P}_t, \mathbf{Q}_t))$ constitutes an i.i.d. sequence with values in $\mathbb{R}^{d \times d} \times \mathbb{R}^d$. A linear SRE is called *polynomial* if there exists an i.i.d. sequence $(\mathbf{e}_t)$ in $\mathbb{R}^{d'}$ such that $\mathbf{P}_t = \mathbf{P}(\mathbf{e}_t)$ and $\mathbf{Q}_t = \mathbf{Q}(\mathbf{e}_t)$, where $\mathbf{P}(\mathbf{x})$ and $\mathbf{Q}(\mathbf{x})$ have entries and coordinates, respectively, which are polynomial functions of the coordinates of $\mathbf{x}$. The existence and uniqueness of a stationary solution to (4.16) has been studied by Brandt [10], Bougerol and Picard [7], Babillot et al. [1] and others. The following set of conditions is sufficient: $E \log^+ \|\mathbf{P}_1\| < \infty$, $E \log^+ |\mathbf{Q_1}| < \infty$, and the top Lyapunov coefficient associated with the operator sequence $(\mathbf{P}_t)$ is strictly negative, that is,

$$(4.17) \qquad \rho = \inf\{t^{-1} E \log \|\mathbf{P}_t \cdots \mathbf{P}_1\| \mid t \geq 1\} < 0.$$

Here $\|\cdot\|$ is the operator norm corresponding to an arbitrary fixed norm $|\cdot|$ in $\mathbb{R}^d$, for example, the Euclidean norm. The following result is a slight generalization of Theorem 4.3 in [21]; see the beginning of the proof below for a comparison.

THEOREM 4.5. *Let $(\mathbf{e}_t)$ be an i.i.d. sequence of random vectors in $\mathbb{R}^{d'}$. Then consider the polynomial linear SRE*

$$(4.18) \qquad \mathbf{Y}_t = \mathbf{P}(\mathbf{e}_t)\mathbf{Y}_{t-1} + \mathbf{Q}(\mathbf{e}_t),$$

*where $\mathbf{P}(\mathbf{e}_t)$ is a random $d \times d$ matrix and $\mathbf{Q}(\mathbf{e}_t)$ a random $\mathbb{R}^d$-valued vector. Suppose:*

1. *$\mathbf{P(0)}$ has spectral radius strictly smaller than 1 and the top Lyapunov coefficient $\rho$ corresponding to $(\mathbf{P}(\mathbf{e}_t))$ is strictly negative.*
2. *There is an $s > 0$ such that*

    $$E\|\mathbf{P}(\mathbf{e}_1)\|^s < \infty \quad and \quad E|\mathbf{Q}(\mathbf{e}_1)|^s < \infty.$$

3. *There is a smooth algebraic variety $V \subset \mathbb{R}^{d'}$ such that $\mathbf{e}_1$ has a density $f$ with respect to Lebesgue measure on $V$. Assume that $\mathbf{0}$ is contained in the closure of the interior of the support $\{f > 0\}$.*

*Then the polynomial linear SRE (4.18) has a unique stationary ergodic solution $(\mathbf{Y}_t)$ which is absolutely regular with geometric rate and consequently strongly mixing with geometric rate.*

REMARK 4.6. As regards the definition of a *smooth algebraic variety*, we first introduce the notion of an *algebraic subset*. An algebraic subset of the $\mathbb{R}^{d'}$ is a set of the form

$$V = \{\mathbf{x} \in \mathbb{R}^{d'} \mid F_1(\mathbf{x}) = \cdots = F_r(\mathbf{x}) = \mathbf{0}\},$$

where $F_1, \ldots, F_r$ are real multivariate polynomials. An *algebraic variety* is an algebraic subset which is not the union of two proper algebraic subsets. An algebraic variety is *smooth* if the Jacobian of $\mathbf{F} = (F_1, \ldots, F_r)^T$ has identical rank everywhere on $V$. Examples of smooth algebraic varieties in $\mathbb{R}^{d'}$ are the hyperplanes of $\mathbb{R}^{d'}$ or $V = \mathbb{R}^{d'}$.

REMARK 4.7. Recall that *absolute regularity* (or *β-mixing*) is a mixing notion which is slightly more restrictive than strong mixing:

$$E\left(\sup_{B \in \sigma(\mathbf{Y}_t, t > k)} |P(B \mid \sigma(\mathbf{Y}_s, s \leq 0)) - P(B)|\right) =: b_k \to 0, \qquad k \to \infty.$$

Indeed, $\beta$-mixing implies strong mixing with the same rate function; see [14] for details on mixing.

PROOF OF THEOREM 4.5. If $E\|\mathbf{P}(\mathbf{e}_1)\|^{\tilde{s}} < 1$ for some $\tilde{s} > 0$, we can immediately apply Theorem 4.3 in [21]. In the general case, we use Mokkadem's result to prove absolute regularity with geometric rate for some subsequence $(\widetilde{\mathbf{Y}}_t) = (\mathbf{Y}_{tm})_{t \in \mathbb{Z}}$, some $m \geq 1$, by observing that $(\widetilde{\mathbf{Y}}_t)$ satisfies the linear SRE (4.19) below. The subsequence argument works because the mixing coefficient $b_k$ is nonincreasing and since $(\mathbf{Y}_t)$ is a Markov process. Then one has the simpler representation

$$b_k = E\left(\sup_{B \in \sigma(\mathbf{Y}_{k+1})} |P(B \mid \sigma(\mathbf{Y}_0)) - P(B)|\right);$$

see, for example, [9].

Since $\rho < 0$, there is an $m \geq 1$ with $E\log\|\mathbf{P}(\mathbf{e}_m)\cdots\mathbf{P}(\mathbf{e}_1)\| < 0$. From the fact that the map $u \mapsto E\|\mathbf{P}(\mathbf{e}_m)\cdots\mathbf{P}(\mathbf{e}_1)\|^u$ has first derivative equal to $E\log\|\mathbf{P}(\mathbf{e}_m)\cdots\mathbf{P}(\mathbf{e}_1)\|$ at $u = 0$, we deduce that there is an $0 < \tilde{s} \leq s$ with $E\|\mathbf{P}(\mathbf{e}_m)\cdots\mathbf{P}(\mathbf{e}_1)\|^{\tilde{s}} < 1$. Then note that $(\widetilde{\mathbf{Y}}_t) = (\mathbf{Y}_{tm})$ obeys a linear SRE:

$$(4.19) \qquad \widetilde{\mathbf{Y}}_t = \widetilde{\mathbf{P}}(\tilde{\mathbf{e}}_t)\widetilde{\mathbf{Y}}_{t-1} + \widetilde{\mathbf{Q}}(\tilde{\mathbf{e}}_t),$$

where

$$\tilde{\mathbf{e}}_t = \begin{pmatrix} \mathbf{e}_{tm} \\ \vdots \\ \mathbf{e}_{(t-1)m+1} \end{pmatrix}$$

and

$$\widetilde{\mathbf{P}}(\tilde{\mathbf{e}}_t) = \mathbf{P}(\mathbf{e}_{tm})\cdots\mathbf{P}(\mathbf{e}_{(t-1)m+1}),$$

$$\widetilde{\mathbf{Q}}(\tilde{\mathbf{e}}_t) = \mathbf{Q}(\mathbf{e}_{tm}) + \sum_{j=1}^{m-1}\left(\prod_{i=1}^{j}\mathbf{P}(\mathbf{e}_{tm+1-i})\right)\mathbf{Q}(\mathbf{e}_{tm-j}).$$

Since both the matrix $\widetilde{\mathbf{P}}(\tilde{\mathbf{e}}_t)$ and the vector $\widetilde{\mathbf{Q}}(\tilde{\mathbf{e}}_t)$ are polynomial functions of the coordinates of $\tilde{\mathbf{e}}_t$ and the sequence $(\tilde{\mathbf{e}}_t)$ is i.i.d., $(\widetilde{\mathbf{Y}}_t)$ obeys a polynomial linear SRE. Observe that $\widetilde{\mathbf{P}}(\mathbf{0}) = (\mathbf{P}(\mathbf{0}))^m$ has spectral radius strictly smaller than 1, that $E\|\widetilde{\mathbf{P}}(\tilde{\mathbf{e}}_1)\|^{\tilde{s}} < 1$ and $E\|\widetilde{\mathbf{Q}}(\tilde{\mathbf{e}}_1)\|^{\tilde{s}} < \infty$ and that $\tilde{\mathbf{e}}_1$ has a density with respect to Lebesgue measure on $V^m$, where $V^m$ is a smooth algebraic variety (see A.14 in [21]). Thus, an application of Theorem 4.3 in [21] yields that $(\widetilde{\mathbf{Y}}_t)$ is absolutely regular with geometric rate. This proves the assertion. $\square$

The following two facts will also be needed.

LEMMA 4.8.  *Let $(\mathbf{P}_t)$ be an i.i.d. sequence of $k \times k$ matrices with $E\|\mathbf{P}_1\|^s < \infty$ for some $s > 0$. Then the associated top Lyapunov coefficient $\rho < 0$ if and only if there exist $c > 0$, $\tilde{s} > 0$ and $\lambda < 1$ so that*

$$(4.20) \qquad E\|\mathbf{P}_t \cdots \mathbf{P}_1\|^{\tilde{s}} \le c\lambda^t, \qquad t \ge 1.$$

PROOF.   For the proof of necessity, observe that there exists $n \ge 1$ such that $E \log \|\mathbf{P}_n \cdots \mathbf{P}_1\| < 0$. From the fact that the map $u \mapsto E\|\mathbf{P}_n \cdots \mathbf{P}_1\|^u$ has first derivative equal to $E \log \|\mathbf{P}_n \cdots \mathbf{P}_1\|$ at $u = 0$, we deduce that there is an $\tilde{s} > 0$ with $E\|\mathbf{P}_n \cdots \mathbf{P}_1\|^{\tilde{s}} = \tilde{\lambda} < 1$. Since the operator norm $\|\cdot\|$ is submultiplicative and the factors in $\mathbf{P}_t \cdots \mathbf{P}_1$ are i.i.d.,

$$E\|\mathbf{P}_t \cdots \mathbf{P}_1\|^{\tilde{s}} \le \tilde{\lambda}^{t/n-1}\left(\max_{\ell=1,\dots,n-1} E\|\mathbf{P}_\ell \cdots \mathbf{P}_1\|^{\tilde{s}}\right) \le c\lambda^t, \qquad t \ge 1,$$

for $c = \tilde{\lambda}^{-1}(\max_{\ell=1,\dots,n-1} E\|\mathbf{P}_\ell \cdots \mathbf{P}_1\|^{\tilde{s}})$ and $\lambda = \tilde{\lambda}^{1/n}$. Regarding the proof of sufficiency, use Jensen's inequality and $\lim_{t \to \infty} t^{-1} E \log \|\mathbf{P}_t \cdots \mathbf{P}_1\| = \rho$ to conclude

$$\rho = \lim_{t \to \infty} \frac{1}{t\tilde{s}} E \log \|\mathbf{P}_t \cdots \mathbf{P}_1\|^{\tilde{s}} \le \limsup_{t \to \infty} \frac{1}{t\tilde{s}} \log E\|\mathbf{P}_t \cdots \mathbf{P}_1\|^{\tilde{s}}$$

$$\le \limsup_{t \to \infty} \frac{1}{t\tilde{s}}(\log c + t \log \lambda) = \frac{\log \lambda}{\tilde{s}} < 0.$$

This completes the proof of the lemma.   □

LEMMA 4.9.  *Suppose that*

$$(4.21) \qquad \mathbf{P}_t = \begin{pmatrix} \mathbf{A}_t & \mathbf{0}_{r \times (k-r)} \\ \mathbf{B}_t & \mathbf{C}_t \end{pmatrix}, \qquad t \in \mathbb{Z},$$

*forms an i.i.d. sequence of $k \times k$ matrices with $E\|\mathbf{P}_1\|^s < \infty$, $s > 0$, where $\mathbf{A}_t \in \mathbb{R}^{r \times r}$, $\mathbf{B}_t \in \mathbb{R}^{(k-r) \times r}$ and $\mathbf{C}_t \in \mathbb{R}^{(k-r) \times (k-r)}$. Then its associated top Lyapunov coefficient $\rho_{\mathbf{P}} < 0$ if and only if the sequences $(\mathbf{A}_t)$ and $(\mathbf{C}_t)$ have top Lyapunov coefficients $\rho_{\mathbf{A}} < 0$ and $\rho_{\mathbf{C}} < 0$.*

PROOF.   For the proof of sufficiency of $\rho_{\mathbf{A}} < 0$ and $\rho_{\mathbf{C}} < 0$ for $\rho_{\mathbf{P}} < 0$, it is by Lemma 4.8 enough to derive a moment inequality of the form (4.20) for $(\mathbf{P}_t)$. By induction we obtain

$$\mathbf{P}_t \cdots \mathbf{P}_1 = \begin{pmatrix} \mathbf{A}_t \cdots \mathbf{A}_1 & \mathbf{0}_{r \times (k-r)} \\ \mathbf{Q}_t & \mathbf{C}_t \cdots \mathbf{C}_1 \end{pmatrix},$$

where

$$\mathbf{Q}_t = \mathbf{B}_t \mathbf{A}_{t-1} \cdots \mathbf{A}_1 + \mathbf{C}_t \mathbf{B}_{t-1} \mathbf{A}_{t-2} \cdots \mathbf{A}_1 + \mathbf{C}_t \mathbf{C}_{t-1} \mathbf{B}_{t-2} \mathbf{A}_{t-3} \cdots \mathbf{A}_1$$

$$+ \cdots + \mathbf{C}_t \cdots \mathbf{C}_3 \mathbf{B}_2 \mathbf{A}_1 + \mathbf{C}_t \cdots \mathbf{C}_2 \mathbf{B}_1.$$

Observe that

(4.22)
$$\max(\|\mathbf{A}_t \cdots \mathbf{A}_1\|, \|\mathbf{C}_t \cdots \mathbf{C}_1\|)$$
$$\leq \|\mathbf{P}_t \cdots \mathbf{P}_1\| \leq \|\mathbf{A}_t \cdots \mathbf{A}_1\| + \|\mathbf{C}_t \cdots \mathbf{C}_1\| + \|\mathbf{Q}_t\|.$$

It is sufficient to show (4.20) for each block in the matrix $\mathbf{P}_t \cdots \mathbf{P}_1$. Because of $\rho_{\mathbf{A}} < 0$, $\rho_{\mathbf{C}} < 0$ and $E\|\mathbf{A}_1\|^s, E\|\mathbf{C}_1\|^s \leq E\|\mathbf{P}_1\|^s < \infty$, Lemma 4.8 already implies moment bounds of the form (4.20) for $(\mathbf{A}_t)$ and $(\mathbf{C}_t)$. Thus, we are left to bound $\|\mathbf{Q}_t\|$. Without loss of generality, we may assume that the constants $\lambda < 1$ and $\tilde{s}, c > 0$ in (4.20) are equal for $(\mathbf{A}_t)$ and $(\mathbf{C}_t)$ and that $\tilde{s} \leq s \leq 1$. From an application of the Minkowski inequality and exploiting the independence of the factors in each summand of $\mathbf{Q}_t$, we obtain the desired relation

$$E\|\mathbf{Q}_t\|^{\tilde{s}} \leq c^2 t E\|\mathbf{B}_1\|^{\tilde{s}} \lambda^{t-1} \leq \tilde{c}\tilde{\lambda}^t,$$

for some $\tilde{\lambda} \in (\lambda, 1)$, $\tilde{c} > 0$. For the proof of necessity, assume $\rho_{\mathbf{P}} < 0$. Then the left-hand side estimates in (4.22) and Lemma 4.8 imply that $\rho_{\mathbf{A}} < 0$ and $\rho_{\mathbf{C}} < 0$. $\square$

We now exploit Theorem 4.5 in order to establish strong mixing with geometric rate of the sequence $(\mathbf{Y}_t) = (\mathbf{G}_t Y_t)$, where $\mathbf{G}_t = h_t'(\boldsymbol{\theta}_0)/\sigma_t^2$ and $Y_t = (Z_t^2 - 1)/2$.

PROPOSITION 4.10. *Let $(X_t)$ be a stationary GARCH$(p,q)$ process with true parameter vector $\boldsymbol{\theta}_0$. Moreover, assume that $Z_1$ has a Lebesgue density $f$, where the closure of the interior of the support $\{f > 0\}$ contains the origin. Then $(\mathbf{Y}_t)$ is absolutely regular with geometric rate.*

PROOF. For the proof of this result, we first embed $(\mathbf{Y}_t)$ in a polynomial linear SRE. Without loss of generality, assume $p, q \geq 3$. Write

$$\widetilde{\mathbf{Y}}_t = \Big(\sigma_{t+1}^2, \ldots, \sigma_{t-q+2}^2, X_t^2, \ldots, X_{t-p+2}^2,$$

$$\frac{\partial h_{t+1}(\boldsymbol{\theta}_0)}{\partial \alpha_0}, \ldots, \frac{\partial h_{t-q+2}(\boldsymbol{\theta}_0)}{\partial \alpha_0}, \ldots, \frac{\partial h_{t+1}(\boldsymbol{\theta}_0)}{\partial \alpha_p}, \ldots, \frac{\partial h_{t-q+2}(\boldsymbol{\theta}_0)}{\partial \alpha_p},$$

$$\frac{\partial h_{t+1}(\boldsymbol{\theta}_0)}{\partial \beta_1}, \ldots, \frac{\partial h_{t-q+2}(\boldsymbol{\theta}_0)}{\partial \beta_1}, \ldots, \frac{\partial h_{t+1}(\boldsymbol{\theta}_0)}{\partial \beta_q}, \ldots, \frac{\partial h_{t-q+2}(\boldsymbol{\theta}_0)}{\partial \beta_q}\Big)^T.$$

Since $Z_t^2 = X_t^2/\sigma_t^2$, we have

$$\sigma(\mathbf{Y}_t, t > k) \subset \sigma(\widetilde{\mathbf{Y}}_t, t > k) \quad \text{and} \quad \sigma(\mathbf{Y}_t, t \leq 0) \subset \sigma(\widetilde{\mathbf{Y}}_t, t \leq 0).$$

Consequently, it is enough to demonstrate absolute regularity with geometric rate of the sequence $(\widetilde{\mathbf{Y}}_t)$. We introduce various matrices. Write $\mathbf{0}_{d_1 \times d_2}$ for the $d_1 \times d_2$ matrix with all entries equal to zero and let $\mathbf{I}_d$ denote the identity matrix of dimension $d$. Then set

$$\mathbf{M}_1(Z_t) = \begin{pmatrix} \boldsymbol{\tau}_t & \beta_q^\circ & \boldsymbol{\alpha}^\circ & \alpha_p^\circ \\ \mathbf{I}_{q-1} & \mathbf{0}_{(q-1)\times 1} & \mathbf{0}_{(q-1)\times(p-2)} & \mathbf{0}_{(q-1)\times 1} \\ \boldsymbol{\xi}_t & \mathbf{0}_{1\times 1} & \mathbf{0}_{1\times(p-2)} & \mathbf{0}_{1\times 1} \\ \mathbf{0}_{(p-2)\times(q-1)} & \mathbf{0}_{(p-2)\times 1} & \mathbf{0}_{(p-2)\times(p-2)} & \mathbf{0}_{(p-2)\times 1} \end{pmatrix},$$

where

$$\boldsymbol{\tau}_t = (\beta_1^\circ + \alpha_1^\circ Z_t^2, \beta_2^\circ, \dots, \beta_{q-1}^\circ) \in \mathbb{R}^{q-1},$$

$$\boldsymbol{\xi}_t = (Z_t^2, 0, \dots, 0) \in \mathbb{R}^{q-1},$$

$$\boldsymbol{\alpha}^\circ = (\alpha_2^\circ, \dots, \alpha_{p-1}^\circ) \in \mathbb{R}^{p-2}.$$

Moreover, define

$$\mathbf{M}_2(Z_t) = \begin{pmatrix} \mathbf{0}_{q\times(p+q-1)} \\ \mathbf{U}_1 \\ \vdots \\ \mathbf{U}_p \end{pmatrix} \quad \text{and} \quad \mathbf{M}_4 = \begin{pmatrix} \mathbf{V}_1 \\ \vdots \\ \mathbf{V}_q \end{pmatrix},$$

where $\mathbf{U}_i \in \mathbb{R}^{q\times(p+q-1)}$ and $\mathbf{V}_j \in \mathbb{R}^{q\times(p+q-1)}$ are given by

$$[\mathbf{U}_1]_{k,\ell} = \delta_{k\ell,11} Z_t^2,$$

$$[\mathbf{U}_i]_{k,\ell} = \delta_{k\ell,1(q+i-1)}, \qquad i \geq 2,$$

$$[\mathbf{V}_j]_{k,\ell} = \delta_{k\ell,1j}.$$

Here $\delta.$ denotes the Kronecker symbol. Also introduce the $q \times q$ matrix

$$\mathbf{C} = \begin{pmatrix} \beta_1^\circ & \cdots & \beta_q^\circ \\ \mathbf{I}_{q-1} & & \mathbf{0}_{(q-1)\times 1} \end{pmatrix},$$

and let

$$\mathbf{M}_3 = \mathrm{diag}(\mathbf{C}, p+1), \qquad \mathbf{M}_5 = \mathrm{diag}(\mathbf{C}, q)$$

be the block diagonal matrices consisting of $p+1$ (or $q$) copies of the block $C$. Finally, we define

$$\mathbf{P}(Z_t) = \begin{pmatrix} \mathbf{M}_1(Z_t) & \mathbf{0}_{(p+q-1)\times(p+1)q} & \mathbf{0}_{(p+q-1)\times q^2} \\ \mathbf{M}_2(Z_t) & \mathbf{M}_3 & \mathbf{0}_{(p+1)q\times q^2} \\ \mathbf{M}_4 & \mathbf{0}_{q^2\times(p+1)q} & \mathbf{M}_5 \end{pmatrix}$$

and $\mathbf{Q} \in \mathbb{R}^{p+q-1+q(p+q+1)}$ by $[\mathbf{Q}]_k = \alpha_0 \delta_{k,1} + \delta_{k,p+q}$. Differentiating both sides of (4.8) at the true parameter $\boldsymbol{\theta} = \boldsymbol{\theta}_0$, we recognize that

$$h'_{t+1}(\boldsymbol{\theta}_0) = (1, X_t^2, \dots, X_{t+1-p}^2, \sigma_t^2, \dots, \sigma_{t+1-q}^2)^T$$
$$+ \beta_1^\circ h'_t(\boldsymbol{\theta}_0) + \cdots + \beta_q^\circ h'_{t+1-q}(\boldsymbol{\theta}_0).$$

From this recursive relationship together with $\sigma_{t+1}^2 = \alpha_0^\circ + \alpha_1^\circ X_t^2 + \cdots + \alpha_p^\circ X_{t+1-p}^2 + \beta_1^\circ \sigma_t^2 + \cdots + \beta_q^\circ \sigma_{t+1-q}^2$, we derive a polynomial linear SRE for $(\widetilde{\mathbf{Y}}_t)$:

$$(4.23) \qquad \widetilde{\mathbf{Y}}_t = \mathbf{P}(Z_t)\widetilde{\mathbf{Y}}_{t-1} + \mathbf{Q}.$$

The proof of Proposition 4.10 follows from the following lemma. $\square$

LEMMA 4.11. *Under the assumptions of Proposition 4.10, the polynomial linear SRE (4.23) has a strictly stationary solution $(\widetilde{\mathbf{Y}}_t)$ which is absolutely regular with geometric rate.*

PROOF. The aim is to show that (4.23) obeys the conditions of Theorem 4.5. Since $EZ_1^2 = 1$, it is immediate that $E\|\mathbf{P}(Z_1)\| < \infty$ since this statement is true for the Frobenius norm and all matrix norms are equivalent. Treat the blocks $\mathbf{M}_1(Z_t)$, $\mathbf{M}_3$ and $\mathbf{M}_4$ separately. Observe that the matrix $\mathbf{M}_1(Z_t)$ appears in the linear SRE for the vector $\mathbf{S}_t = (\sigma_{t+1}^2, \ldots, \sigma_{t-q+2}^2, X_t^2, \ldots, X_{t-p+2}^2)^T$, namely,

$$\mathbf{S}_t = \mathbf{M}_1(Z_t)\mathbf{S}_{t-1} + (\alpha_0^\circ, 0, \ldots, 0)^T.$$

Theorem 1.3 of [6] says that (1.1) admits a unique stationary solution if and only if $(\mathbf{M}_1(Z_t))$ has strictly negative top Lyapunov coefficient; consequently, $\rho_{\mathbf{M}_1} < 0$. Moreover, arguing by recursion on $p$ and expanding the determinant with respect to the last column, it is easily verified that $\mathbf{M}_1(0)$ has characteristic polynomial

$$\det(\lambda \mathbf{I}_{p+q-1} - \mathbf{M}_1(0)) = \lambda^{p+q-1}\left(1 - \sum_{i=1}^q \beta_i^\circ \lambda^{-i}\right).$$

Since (4.1) holds for a stationary GARCH$(p,q)$ process, by the triangle inequality

$$\left|1 - \sum_{i=1}^q \beta_i^\circ \lambda^{-i}\right| \geq 1 - \sum_{i=1}^q \beta_i^\circ \lambda^{-i} \geq 1 - \sum_{i=1}^q \beta_i^\circ > 0$$

if $|\lambda| \geq 1$ and, hence, $\mathbf{M}_1(0)$ has spectral radius $< 1$. Observe that the building block $\mathbf{C}$ has characteristic polynomial

$$\det(\lambda \mathbf{I}_q - \mathbf{C}) = \lambda^q\left(1 - \sum_{i=1}^q \beta_i^\circ \lambda^{-i}\right),$$

showing that its spectral radius is strictly smaller than 1 (use the same argument as before). Thus, the *deterministic* matrices $\mathbf{M}_3$ and $\mathbf{M}_5$ have spectral radius $< 1$, which also implies that their associated top Lyapunov coefficients are strictly negative. Combining these results, we deduce that $\mathbf{P}(0)$ has spectral radius $< 1$ and conclude by twice applying Lemma 4.9 that $(\mathbf{P}(Z_t))$ has strictly negative top Lyapunov coefficient. Hence, by Theorem 4.5 the stationary sequence $(\widetilde{\mathbf{Y}}_t)$ is absolutely regular with geometric rate. □

REMARK 4.12. Since $(X_t^2, \sigma_t^2)$ is a subvector of $\widetilde{\mathbf{Y}}_t$, stationary GARCH$(p, q)$ processes are absolutely regular with geometric rate; this result has previously been established by Boussama [8].

**Acknowledgments.** Daniel Straumann would like to thank RiskLab Zürich for the opportunity to conduct this research. In particular, he is grateful to Paul Embrechts and Alexander McNeil for stimulating discussions. Both authors would like to thank the two referees for various remarks which led to a better presentation of the paper.

LABORATORY OF ACTUARIAL MATHEMATICS
DEPARTMENT OF APPLIED MATHEMATICS
  AND STATISTICS
UNIVERSITY OF COPENHAGEN
UNIVERSITETSPARKEN 5
DK-2100 COPENHAGEN Ø
DENMARK
AND
MAPHYSTO
THE DANISH RESEARCH NETWORK
  FOR MATHEMATICAL PHYSICS
  AND STOCHASTICS
E-MAIL: mikosch@math.ku.dk

RISKLAB
DEPARTMENT OF MATHEMATICS
ETH ZÜRICH
ETH ZENTRUM
CH-8092 ZÜRICH
SWITZERLAND
E-MAIL: straumann@math.ethz.ch